\documentclass[11pt]{article}
\usepackage{amsmath,amsthm,amssymb,ifthen,graphicx}
\usepackage[sort]{natbib}
\usepackage{pstricks,pst-node}
\usepackage{subfigure}
\newcounter{isamac} 
\ifthenelse{\value{isamac}=1}{
\input BoxedEPS          %
\SetTexturesEPSFSpecial  
\HideDisplacementBoxes   %
}{}

\newtheorem{theorem}{Theorem}
\newtheorem{algorithm}[theorem]{Algorithm}

\newtheorem{remark}[theorem]{Remark}
\newtheorem{example}[theorem]{Example}

\newcommand{\RRR}{ \mathbb{R}} 
 
\newcommand{\ND}{ \mathcal{N}} 
\newcommand{\I}[2]{ {#1}{\ind}{#2}} 

\newcommand{\tr}{\mathrm{tr}}
\newcommand{\cvec}{\mathrm{vec}}
\newcommand{\bi}{\leftrightarrow}

\newcommand{\Ew}{\text{\rm E}}       
\newcommand{\spo}{\text{\rm sp}}      
\newcommand{\nsp}{\text{\rm nsp}}      

\newcommand{\indm}[2]{\ensuremath{{\mathfrak I}_{\kern-1pt\scriptstyle#1}({\mathcal
#2})}} 

\newcommand{\ind}{\mbox{$\perp \kern-5.5pt \perp$}}

\newcommand{\uned}{\hbox{\kern3pt\raise2.5pt\vbox{\hrule
width9pt height 0.3pt}\kern3pt}}

\newcommand{\dashed}{\hbox{\kern3.05pt\raise2.5pt\vbox{\hrule
width1.7pt height 0.3pt}\kern1.8pt\raise2.5pt\vbox{\hrule
width1.7pt height 0.3pt}\kern1.8pt\raise2.5pt\vbox{\hrule
width1.7pt height 0.3pt}\kern1.8pt\raise2.5pt\vbox{\hrule
width1.7pt height 0.3pt}\kern3.05pt}}

\newcommand{\pedg}[2]{\ensuremath{{\kern0.5pt
\scriptstyle{\ifthenelse{\equal{\head}{#1}}{\lhead\kern0.5pt}{#1\kern0.5pt}}\joinrel\relbar
\negthinspace\relbar\joinrel{\kern0.5pt #2}\kern0.5pt}}}

\newcommand{\pdots}{\hbox{\kern2.5pt\raise1.5pt\hbox{\ensuremath{\ldots}}\ke
rn2.5pt}}  

\newcommand{\h}[1]{\hat{#1}}

\numberwithin{equation}{section}

\begin{document}
 
\renewcommand{\labelenumi}{(\roman{enumi})}
\renewcommand{\thefootnote}{\fnsymbol{footnote}}

\title{Estimation of a Covariance Matrix with 
  Zeros\protect\footnotemark} \author{
  Sanjay Chaudhuri\\
  \em University of Washington\\
  \and
  Mathias Drton\thanks{Corresponding author: {\tt
    drton@galton.uchicago.edu}} \\
  \em University of Chicago
  \and Thomas S.\ Richardson\\
  \em University of Washington} \maketitle
\footnotetext{\protect\footnotemark
  Material in \S \ref{newmethod} appeared as {\em A new algorithm
    for maximum likelihood estimation in Gaussian graphical models for
    marginal independence\/} in the Proceedings of the Conference on
    Uncertainty in Artificial Intelligence, 2003.}


\begin{abstract}
  We consider estimation of the covariance matrix of a multivariate random
  vector under the constraint that certain covariances are zero.  We first
  present an algorithm, which we call Iterative 
  Conditional Fitting, for computing the maximum likelihood estimator of
  the constrained covariance matrix, under the assumption of multivariate
  normality. In contrast to previous approaches, this algorithm has
  guaranteed convergence properties.  Dropping the assumption of
  multivariate normality, we show how to estimate the covariance matrix in
  an empirical likelihood approach.
  These approaches are then compared via simulation and on an example of
  gene expression.\\[0.1cm]
  {\scriptsize {\em Some key words}: Covariance graphs; Empirical
    likelihood; Graphical models; Marginal independence; Maximum
    likelihood estimation; Multivariate normal distribution }
\end{abstract}


\section{Introduction}
\begin{figure}[t]
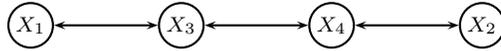

    \vspace{3cm}\hspace{3.5cm}
    \scriptsize
    \psset{linewidth=0.8pt}          
    \newlength{\MyLength}
    \settowidth{\MyLength}{$X$}
    \newcommand{\myNode}[2]{\circlenode{#1}{\makebox[\MyLength]{#2}}}
    \rput(6, 2){\myNode{2}{$X_2$}} 
    \rput(0, 2){\myNode{1}{$X_1$}}
    \rput(4, 2){\myNode{4}{$X_4$}}
    \rput(2, 2){\myNode{3}{$X_3$}} 
    \ncline{<->}{1}{3}
    \ncline{<->}{3}{4}
    \ncline{<->}{2}{4} 
    \vspace{-1.5cm}
\caption{\label{fig:ex.bi} The covariance graph for the matrix in
  (\ref{eq:sigmaex}).}
\end{figure}

In this paper we consider estimation of the covariance matrix of a random
vector, subject to certain entries being set to zero. Such restrictions
appear, for example, in recent work by \cite{grzebyk:2004} and
\cite{mao:2004}. Suppose we have a random vector
$X=(X_1,X_2,X_3,X_4)'\in \RRR^4$ whose covariance matrix $\Sigma$ exhibits
the zero pattern
\begin{equation}
  \label{eq:sigmaex}
  \Sigma = 
  \begin{pmatrix}
    \sigma_{11} & 0 & \sigma_{13} & 0\\
    0 & \sigma_{22} & 0 & \sigma_{24}\\
    \sigma_{13} & 0 & \sigma_{33} & \sigma_{34}\\
    0 & \sigma_{24} & \sigma_{34} & \sigma_{44}
  \end{pmatrix}\in\RRR^{4\times 4}.
\end{equation}
It is often helpful to visualize the pattern of zeros by a so-called
covariance graph, especially for larger covariance matrices
\citep{coxwerm:lindep,coxwerm:book}. A covariance graph has one vertex
for each one of the random variables in the random vector. In the
above example, the vertex set is $V=\{1,2,3,4\}$, where the random
variable $X_i$ is identified with its index $i$.  Next, each pair of
vertices $(i,j)\in V\times V$, $i\not= j$, is connected by an edge
unless $\sigma_{ij}=0$.  Assuming that the covariance matrix in
(\ref{eq:sigmaex}) has no zeros other than those indicated explicitly,
its covariance graph is given in Figure \ref{fig:ex.bi}.  Here we use
bi-directed edges in keeping with the path diagram notation used by
\citet{wright:cnc}; other authors have used dashed edges; see
\citet{coxwerm:lindep,coxwerm:book}.

We define a {\it covariance graph model} as the set of joint
distributions in which the associated zero restrictions hold in the
covariance matrix. In the absence of an assumption of normality, the
model does not have a Markov interpretation.

The {\it Gaussian covariance graph model} is the family of
all multivariate normal distributions $\ND(\mu,\Sigma)$ such that
$\sigma_{ij}=0$ whenever $i\not= j$ and $i\not\bi j$.  Clearly,
$\sigma_{ij}=0$ if and only if $X_i$ and $X_j$ are marginally
independent; in symbols $\I{X_i}{X_j}$.  Hence a Gaussian covariance
graph model is a graphical model based on marginal independence
in contrast with graphical models based on undirected graphs (Markov
random fields), directed acyclic graphs (DAGs, Bayesian networks), or
chain graphs, where the absence of an edge between two vertices
generally indicates some conditional independence between the
associated variables \citep{edwards:2000,lau:bk,whittaker:bk}.

%

Maximum likelihood (ML) estimation in Gaussian covariance graph models
is not well developed: the conceptual simplicity of these models
belies the fact that, in contrast to undirected graph models, they
form curved exponential families.  For instance, the graphical
modelling software MIM \citep[\S7.4]{edwards:2000} permits fitting of
such models only by a heuristic ``dual likelihood'' method due to
\cite{kauermann:dual}.  There is, however, an algorithm due to
\citet{anderson:1969,anderson:1970,anderson:1973} that can be used to
compute the ML estimate in models defined by linear hypotheses on
covariance matrices, hence also in covariance graph models.  However,
it is unclear when this algorithm converges and when its limit points
are positive semi-definite matrices.  
Such issues become more pressing when mis-specified models are fitted, as
will be the case in a specification search.
In this paper, we introduce a
new algorithm for ML estimation in covariance graph models, called
Iterative Conditional Fitting (ICF), which does not suffer from the
same problems as Anderson's algorithm.

For situations in which multivariate normality does not hold,
estimates may still be obtained via procedures based on normality such
as ICF and dual estimation but the behaviour of these methods is then
unclear.  As an alternative, we present an approach, based on
empirical likelihood \citep{owen:2001}, which provides consistent
estimates even without normality.  We compare the different
estimation methods on real and simulated data.


\section{Covariance graph models}\label{covgraphmodel}

\subsection{Non-parametric model}
\label{sec:nonpar}

Suppose that we observe a random vector $Y_V= (Y_i\mid i\in
V)'\in\RRR^V$, indexed by $V$, and with joint distribution $P_V$. Let
$\Sigma(P_V)=(\sigma_{ij})\in\RRR^{V\times V}$ be the unknown
covariance matrix. Let $G=(V,E)$ be a graph with the variable set $V$
as vertex set and the edge set $E\subseteq V\times V\setminus
\{(i,i)\mid i\in V\}$ consisting exclusively of bi-directed edges
$(i,j),(j,i)\in E$, denoted by $i\bi j$.  Let $\mathbf{P}(V)$ be the
cone of positive definite $V\times V$ matrices and let $\mathbf{P}(G)$
be the cone of all matrices $\Sigma\in\mathbf{P}(V)$ which fulfill the
linear restrictions
\begin{equation}\label{covrestrictions}
i\not\bi j\quad\Longrightarrow\quad\sigma_{ij}=0.
\end{equation}
The {\em covariance graph model\/} $\mathbf{M}(G)$ 
associated with the bi-directed
graph $G$ is simply the family of joint distributions
\begin{equation}\label{defcovgraph}
\mathbf{M}(G)= \big\{ P_V \mid  \Sigma(P_V) \in
\mathbf{P}(G)\big\}. 
\end{equation}

We consider estimation of the unknown parameter $\Sigma=\Sigma(P_V)$
based on a sample of observations $Y^{(k)}_V\in\RRR^V$, $k\in
N=\{1,\ldots, n\}$, that are i.i.d.~according to $P_V\in
\mathbf{M}(G)$.  The set $N$ can be interpreted as indexing the
subjects on which we observe the variables in $V$. We group the
vectors in the sample as columns in the $V\times N$ random matrix $Y$
so that
\begin{equation}
{\rm Var}(Y) = \Sigma\otimes I_N.
\end{equation}
Here, $I_N$ is the $N\times N$ identity matrix and $\otimes$ is the
Kronecker product.  Thus the $i$-th row $Y_{i}\in\RRR^N$ of the
matrix $Y$ contains the i.i.d.~observations for variable $i\in V$ on
all the subjects in $N$ and the $k$-th column $Y^{(k)}_V$ holds all
the observations made on subject $k\in N$. Finally, the sample size is
$n=|N|$ and the number of variables is $p=|V|$.



\subsection{Gaussian model}

We define a {\it Gaussian covariance graph model} as the multivariate normal
submodel
\begin{equation}\label{defgausscovgraph}
\mathbf{N}(G)= \big( \ND_V(\mu,\Sigma) \mid  \Sigma\in
\mathbf{P}(G)\big) \subset \mathbf{M}(G). 
\end{equation}
The {\em log-likelihood function\/} $\ell$ of the covariance graph model
$\mathbf{N}(G)$ is a function from $\RRR^{V} \times \mathbf{P}(G)$ to $\RRR$ and
can be expressed as 
\begin{equation}\label{loglikelihood}
\ell(\mu,\Sigma)
=-\frac{np}{2}\log(2\pi)-\frac{n}{2}\log
|\Sigma|-\frac{n}{2}\tr(\Sigma^{-1}\tilde S),
\end{equation}
see e.g.~\citet[\S3.1]{edwards:2000}. 
Here $\tilde S$ is 
\begin{equation}
\tilde S=\frac{1}{n} (Y- \mu \otimes 1_N) (Y- \mu \otimes
1_N)'\in\RRR^{V\times V},
\end{equation}
where $1_N=(1,\dots,1)'\in \RRR^N$.  For any given value of $\Sigma$,
(\ref{loglikelihood}) is maximized by setting $\mu= \bar Y\in\RRR^V$,
i.e., the vector of the row means of $Y$.  Hence, the profile
log-likelihood for $\Sigma$, $\ell (\Sigma)$ is obtained by replacing
$\tilde S$ with
\begin{equation}\label{empcov}
S=\frac{1}{n} (Y- \bar Y \otimes 1_N) (Y-\bar Y \otimes
1_N)'\in\RRR^{V\times V},
\end{equation}
in (\ref{loglikelihood}). Working with the profile likelihood
corresponds to fitting the submodel of $\mathbf{N}(G)$ in which $\mu=0$ and
adjusting the sample size to $n-1$.


If $S$ is positive definite, which will occur with probability 1 if
$n\ge p+1$ \citep{eaton:1973}, then the global maximum of
$\ell(\Sigma)$ over $\mathbf{P}(G)$, i.e., the ML estimator of
$\Sigma$, exists. In general, the condition $n\ge p+1$ is not
necessary for almost sure existence of the ML estimator but we are not
aware of any results in the literature which provide a necessary and
sufficient condition \citep[compare][]{buhl:1993}. In the sequel we
will assume $S$ to be positive definite.  Note that since the model
$\mathbf{N}(G)$ is a curved, but not necessarily regular, exponential
family, the likelihood function may, and in fact can, have multiple
local maxima \citep{drton:2004,drton:2004b}.



Let
\begin{equation}
  \label{eq:free}
  F =\{ (i,i) \mid i\in V\}\cup 
  \{ (i,j)\in V^2 \mid i<j \wedge i\bi j \} 
\end{equation}
be the pairs of vertices indexing unrestricted elements in the matrix
$\Sigma\in\mathbf{P}(G)$. The cardinality of $F$ is equal to the number of
vertices plus the number of edges in the graph $G$.
The unrestricted elements of $\Sigma$ form the vector
\begin{equation}
  \label{eq:sigvec}
  \sigma = (\sigma_{ij}\mid (i,j)\in F) \in \RRR^{F}.
\end{equation}
In order to write derivatives of the log-likelihood function in compact
form we introduce the matrix $Q$ with entries in $\{0,1\}$
that satisfies $\cvec(\Sigma)=Q\sigma$, where $\cvec$ is the operator of
column-wise matrix vectorization.  The columns of $Q$ that are
associated with a variance $\sigma_{ii}$ contain exactly one entry equal to
one, whereas a column of $Q$ that is associated with a covariance
$\sigma_{ij}$, $i\not= j$, $i\bi j$, contains exactly two entries equal to
one. If the graph $G$ is complete, 
i.e., all possible edges are present in $G$, then $Q$ is the duplication
matrix described in \citet[p.352]{harville:1997}.  

The first derivative of the
log-likelihood function, that is, the {\em score function\/} can then be
written as 
\begin{equation}
  \label{eq:score}
  \frac{\partial\ell(\Sigma)}{\partial\sigma} =
  \frac{n}{2} Q' \big[ \cvec(\Sigma^{-1}S\Sigma^{-1}) - \cvec(\Sigma^{-1})
  \big],
\end{equation}
see \citet[\S 15]{harville:1997} for details on the necessary matrix
differential calculus. It follows that the {\em likelihood equations\/}
$\partial\ell(\Sigma)/\partial\sigma=0$  are
\begin{equation}\label{likequs}
  (\Sigma^{-1})_{ij} =(\Sigma^{-1}S\Sigma^{-1})_{ij},  \quad (i,j)\in F;
\end{equation}
compare also \citet[\S 2.1.1]{anderson:1985}. 
The full matrix $\Sigma$ is determined by $\sigma_{ij}=0$ for
$(i,j)\not\in F$, that is for $i\not= j$ and $i\not\bi j$.


The second derivative of $\ell(\Sigma)$ can be computed using 
results from \citet[\S 15.9]{harville:1997}, and we find that the Hessian
matrix equals
\begin{equation}
  \label{eq:2deriv}
  \frac{\partial^2\ell(\Sigma)}{\partial\sigma^2} =
  \frac{n}{2} Q' \big\{ [\Sigma^{-1}\otimes\Sigma^{-1}] -
  [(\Sigma^{-1}S\Sigma^{-1})\otimes \Sigma^{-1}] -
  [\Sigma^{-1}\otimes(\Sigma^{-1}S\Sigma^{-1})] \big\} Q.
\end{equation}
Its negated expectation under $\ND_V(0,\Sigma)$, 
the {\em Fisher-information\/}, equals
\begin{equation}
  \label{eq:fisherinfo}
  -\Ew\left[\frac{\partial^2\ell(\Sigma)}{\partial\sigma^2}\right] =
  \frac{n}{2} Q' (\Sigma^{-1}\otimes\Sigma^{-1}) Q
\end{equation}
and can be used for normal approximation to the distribution of roots
of the likelihood equations.  Sections \S\S
\ref{sec:alter}-\ref{sec:icfsur} focus on the computations of such
roots.

\section{Existing estimation methods for Gaussian covariance graphs}
\label{sec:alter}

We are aware of only one specialized algorithm for ML estimation applicable to
covariance graph models.  This algorithm is due to \cite{anderson:1973} and 
will fit any Gaussian model obtained from a linear hypothesis
on the covariance matrix \citep{anderson:1969,anderson:1970}. In this
section, we describe the incarnation of {\em Anderson's algorithm\/} that fits
covariance graph models. We also review a {\em dual estimation\/} method due  
to \cite{kauermann:dual}, which produces estimates that are unique and asymptotically
efficient though, in general,
not solutions to the likelihood equations.
Note that \cite{wermcoxmarc:04} have recently proposed moment based
estimators in the special case where the graph is a chain, or
equivalently, the covariance matrix is tri-diagonal under a suitable
ordering.

\subsection{Anderson's algorithm for ML estimation}\label{andmethod}

Each iteration of Anderson's algorithm solves a system of linear
equations built from the current estimate of $\Sigma$. In the case of
covariance graphs, the linear equations are solved for the vector
$\sigma$ of unrestricted elements in $\Sigma$, compare
(\ref{eq:sigvec}), and can be specified as follows. Let
$\sigma^{ij}=(\Sigma^{-1})_{ij}$ and $A= A_\Sigma$ be the $F\times F$
matrix with entries
\begin{equation}
  \label{eq:A}
  A_{(ij,k\ell)} = \left\{ 
    \begin{array}{lll}
      \sigma^{ik}\sigma^{jk} & \; \textrm{if} \; & k=\ell,\\
      \sigma^{ik}\sigma^{j\ell}+\sigma^{jk}\sigma^{i\ell}  & \;
      \textrm{if} \; & k\not=\ell. 
    \end{array}
    \right.
\end{equation}
Here $(i,j)$ and $(k,\ell)$ are elements of $F$.
Furthermore, let $b=b_\Sigma$ be the $F\times 1$ vector with components
\begin{equation}
b_{ij}=(\Sigma^{-1}S\Sigma^{-1})_{ij},\quad (i,j)\in F.
\end{equation}
From \cite{anderson:1973}, it follows that $\Sigma\in\mathbf{P}(G)$ solves
$A_\Sigma\sigma=b_\Sigma$ if and only if $\Sigma$ solves the likelihood
equations \eqref{likequs}. 

This motivates the following iterative scheme.  Start with some
$\Sigma^{(0)}\in\mathbf{P}(G)$. Iteratively update the current
estimate $\Sigma^{(r)}$ to $\Sigma^{(r+1)}$ determined by the linear
equations 
\begin{equation}
A_{\Sigma^{(r)}}\;\sigma^{(r+1)}=b_{\Sigma^{(r)}}.
\end{equation}
A fixed point of this algorithm solves the likelihood
equations \eqref{likequs}. As starting value, Anderson suggests the
identity matrix, i.e., $\Sigma^{(0)}=I_V$. In the first step, his algorithm
constructs the empirical estimate 
$\Sigma^{(1)}$ with $\sigma^{(1)}_{ij}=S_{ij}$, $(i,j)\in F$. However,
neither $\Sigma^{(1)}$ nor any subsequent estimate of $\Sigma$ has to be
positive {(semi-)} definite and thus may not be a valid
covariance matrix. Moreover, at any given stage, the likelihood function may
decrease, and convergence of Anderson's algorithm cannot be guaranteed.

\subsection{Kauermann's dual estimation}
\label{sec:dual}

Dual estimation is based on the maximization of a dual likelihood 
function, which is motivated by interchanging the role of the
parameter matrix $\Sigma$ and the empirical covariance matrix $S$ in \eqref{loglikelihood}
\citep[\S4]{kauermann:dual}.  Procedurally, dual estimation, also called
minimizing the discriminant information,
amounts to finding the
matrix $\hat\Sigma_\mathrm{dual}\in\mathbf{P}(G)$ that solves the
equations
\begin{equation}
  \label{eq:dualeqn}
  (\hat\Sigma_\mathrm{dual}^{-1})_{ij}  = (S^{-1})_{ij}, \quad \forall
  (i,j)\in F, 
\end{equation}
while satisfying that $(\hat\Sigma_\mathrm{dual})_{ij} = 0$ for all
$(i,j)\not\in F$.  Contrary to (\ref{likequs}), the equation system
(\ref{eq:dualeqn}) always has a unique solution that can be found by the
iterative proportional fitting algorithm; see also
\citet[\S7.4]{edwards:2000}.  In particular, if the covariance graph is
decomposable, then iterative proportional fitting will terminate in
finitely many steps, and the dual estimator $\hat\Sigma_\mathrm{dual}$
is available in closed form.

\section{Iterative conditional fitting for Gaussian covariance
graphs}\label{newmethod}

In this section, we present the new {\em Iterative Conditional Fitting\/}
(ICF) algorithm for ML estimation, which is guaranteed to produce positive
definite roots of the likelihood equations of covariance graph models. We
begin by explaining the idea of iteratively fitting conditional
distributions that stands behind ICF, and then show how the algorithm can
be implemented using simple least squares computations.

\subsection{The idea of iterative conditional fitting}

Starting with some initial estimate of the joint distribution, the
idea of ICF is to repeatedly iterate through all vertices $i\in 
V$, and
\begin{enumerate}
\item\label{step1} Fix the marginal distribution for the variables
  different from 
  $i$, i.e., the variables $-i= V\setminus \{i\}$;
\item\label{step2} Estimate, by maximum likelihood, the conditional
  distribution of 
  variable $i$ given the variables $-i$ under the constraints implied by
  the covariance graph  model $\mathbf{N}(G)$;
\item\label{step3} Find a new estimate of the joint distribution by
  multiplying together the fixed marginal and the estimated conditional
  distribution. 
\end{enumerate}
Since we fix the marginal distribution of variables $-i$ in
the update for variable $i$, all marginal independences
amongst the variables $-i$ still hold true after the update.
Therefore, only the marginal independences involving variable $i$ lead
to constraints for the estimation in step (ii). 

In order to make the idea more precise, let $\Sigma_{A,B}$ 
denote the $A\times B$ submatrix of $\Sigma$ and $Y_A$ denote the
$A\times N$ submatrix of $Y$, where $A,B\subseteq V$. 
Clearly, 
\[
Y_{-i}\sim \ND_{-i\times N}(0,\Sigma_{-i,-i}\otimes I_N).
\]
Hence, step (i) simply fixes the value of
$\Sigma_{-i,-i}$, i.e., everything but the $i$-th row and column of
$\Sigma$.  As $\Sigma_{-i,-i}$ remains unchanged in the $i$-th update many
of the zero constraints imposed on the covariance matrix trivially hold
true also after the update.

The conditional distribution of $Y_i$ given $Y_{-i}$ is the normal
distribution 
\begin{equation}\label{conddist}
(Y_i\mid Y_{-i})\sim 
\ND_{\{i\}\times N}(B_{i}Y_{-i},\lambda_{i}I_N),
\end{equation}
where 
\begin{equation}\label{defB}
B_{i}=\Sigma_{i, -i}(\Sigma_{-i, -i})^{-1}\in\RRR^{\{i\}\times -i}
\end{equation}
is the $\{i\}\times -i$ matrix of regression coefficients, and
\begin{equation}\label{deflambda}
\lambda_{i}= \sigma_{ii}-\Sigma_{i,
-i}(\Sigma_{-i,-i})^{-1}\Sigma_{-i, i}\in (0,\infty)
\end{equation}
is the conditional variance. 
If the
graph $G$ was the complete graph $\bar G$ 
in which an edge joins any pair of vertices then the mapping
\begin{equation}\label{biject1}
\begin{split}
\mathbf{P}(\bar G)=\mathbf{P}(V) &\to (0,\infty) \times
\RRR^{\{i\}\times -i}\times \mathbf{P}_{-i}(\bar G),\\
\Sigma &\mapsto (\lambda_i, B_{i}, \Sigma_{-i,-i})
\end{split}
\end{equation}
would be bijective and the regression in \eqref{conddist} a standard
least squares regression. Here, $\mathbf{P}_A(G)$ is the set of all
$A\times A$ submatrices of matrices in $\mathbf{P}(G)$, $A\subseteq
V$.  For a general graph $G$, (\ref{biject1}) is no longer bijective and
(\ref{conddist}) is not a standard regression because we need to
respect the restriction $\Sigma\in\mathbf{P}(G)$, i.e., the
restrictions $\sigma_{ij}=0$ if $j\in -i$, $j\not\bi i$.  However,
this can be circumvented using synthetic {\em
  pseudo-variables\/} that are computed from the data $Y_{-i}$ and the
fixed matrix $\Sigma_{-i,-i}$.


\subsection{Pseudo-variable regressions}
\label{sec:pseudoreg}

Instead of working with the regressions coefficients $B_i$, we exploit
the fact that $B_i$ equals $\Sigma_{i,-i}$ multiplied by the inverse
of the fixed submatrix $\Sigma_{-i,-i}$.  Let $\spo(i)=\{j\mid i\bi
j\}$ be the set of {\em spouses\/} of $i\in V$ and let
$\nsp(i)=V\setminus(\spo(i)\cup\{i\})$ be the set of {\em
  non-spouses\/}, yielding the partition
$V=\{i\}\cup \spo(i)\cup \nsp(i)$. Then the conditional
expectation of $(Y_i\mid Y_{-i})$ can be written as
\begin{equation}
  \label{eq:1}
  \Ew[Y_i\mid Y_{-i}] = \Sigma_{i,-i}\big[ (\Sigma_{-i,-i})^{-1}\,Y_{-i}\big]
  =  \Sigma_{i,\spo(i)}Z^{(i)}_{\spo(i)} = 
  \sum_{j\in\spo(i)} \sigma_{ij}  Z^{(i)}_j,
\end{equation}
where the {\em pseudo-variable\/} $Z^{(i)}_j$ is equal to the $j$-th row in
\begin{equation}
\label{eq:2}
 Z^{(i)}_{\spo(i)} = [(\Sigma_{-i,-i})^{-1}]_{\spo(i),-i}\,Y_{-i}
 \in\RRR^{\spo(i)\times N}.
\end{equation}
In (\ref{eq:1}), we exploit that $\sigma_{ij}=0$ if $j\in \nsp(i)$.
From (\ref{eq:1}), we obtain 
\begin{equation}
  \label{eq:conddistpseud}
  (Y_i\mid Y_{-i})\sim 
  \ND_{\{i\}\times N}\Big(
   \sum_{j\in\spo(i)} \sigma_{ij}  Z^{(i)}_j 
  ,\lambda_{i}I_N\Big).
\end{equation}
Let $\mathbf{P}_{-i}(G)$ be the set of $-i\times -i$ submatrices of
the matrices in $\mathbf{P}(G)$. Then the mapping 
\begin{equation}\label{biject2}
\begin{split}
\mathbf{P}(G) &\to (0,\infty) \times
\RRR^{\{i\}\times \spo(i)}\times \mathbf{P}_{-i}(G)\\
\Sigma &\mapsto (\lambda_i, \Sigma_{i,\spo(i)}, \Sigma_{-i,-i})
\end{split}
\end{equation}
is a bijection, which implies that the parameters $\sigma_{ij}$, $j\in
\spo(i)$, and $\lambda_i$ are variation independent
in (\ref{eq:conddistpseud}).  Therefore, if  
$\Sigma_{-i,-i}$ is fixed to equal some
given matrix in $\mathbf{P}_{-i}(G)$, then (\ref{eq:conddistpseud})
constitutes a standard normal regression model whose parameters
$\sigma_{ij}$, $j\in \spo(i)$, and $\lambda_i$ can be estimated by
the usual least squares formula.  The estimate of $\lambda_i$ yields
an estimate of $\sigma_{ii}$ by solving (\ref{deflambda}) for
$\sigma_{ii}$.  Thus, we obtain the ML estimator of the $i$-th row
and column of $\Sigma$ when $\Sigma_{-i,-i}$ is fixed.
In particular, after updating the $i$-th row and column we are still left with a matrix
$\Sigma\in\mathbf{P}(G)$. 

\subsection{The iterative conditional fitting algorithm}\label{algorithm}

Let $\hat\Sigma^{(r)}$ be the estimate of $\Sigma$ after the $r$-th
iteration and $\hat\Sigma^{(r,i)}$  the estimate of $\Sigma$ 
after the $i$-th update step of the $r$-th iteration in ICF, i.e., after
estimating $(Y_i\mid Y_{-i})$.
\begin{algorithm}
  \label{alg:icf}
  The ICF algorithm can be implemented as:
\begin{enumerate}
\item[1.] {\em (Initialization)} Set the iteration counter $r=0$, and
  choose a starting value $\hat \Sigma^{(0)}\in\mathbf{P}(G)$,
  e.g.~the identity matrix $\hat \Sigma^{(0)}=I_V$.
\item[2.] {\em (Updates)} Order the variables as $V=\{1,\ldots,p\}$,
  set $\hat \Sigma^{(r,0)}=\hat \Sigma^{(r)}$, and repeat the
  following steps for all $i=1,\ldots,p$ :
\begin{enumerate}
\item[(i)] Let $\hat \Sigma^{(r,i)}_{-i,-i}=\hat
  \Sigma^{(r,i-1)}_{-i,-i}$ and calculate from this submatrix the
  pseudo-variables $Z^{(i)}_{\spo(i)}$ according to (\ref{eq:2}).
\item[(ii)] Compute the ML estimators
\begin{equation}\label{ls}
\begin{split}
  \mspace{-10mu}\hat \Sigma^{(r,i)}_{i,\spo(i)}&= Y_i
  \,(Z^{(i)}_{\spo(i)})'\big[ Z^{(i)}_{\spo(i)}(Z^{(i)}_{\spo(i)})'
  \big]^{-1}
  ,\\
  \hat \lambda_i &= \frac{1}{n}(Y_i-\hat
  \Sigma^{(r,i)}_{i,\spo(i)}Z^{(i)}_{\spo(i)})(Y_i-\hat
  \Sigma^{(r,i)}_{i,\spo(i)}Z^{(i)}_{\spo(i)})'.\\
\end{split}
\end{equation}
for the linear regression (\ref{eq:conddistpseud}). The existence of
the matrix inverse follows from the assumed non-singularity of the
sample covariance matrix $S$.
\item[(iii)] Complete $\hat\Sigma^{(r,i)}$ by
setting 
\begin{equation}
  \hat \sigma^{(r,i)}_{ii}=\hat
  \lambda_i+\hat\Sigma^{(r,i)}_{i,\spo(i)}
  \big[(\hat\Sigma^{(r,i)}_{-i,-i})^{-1}\big]_{\spo(i),\spo(i)}
  \hat\Sigma^{(r,i)}_{\spo(i),i}\;;  
\end{equation}
compare (\ref{deflambda}). 
\end{enumerate}
\item[3.] {\em (Repeat)} Set $\hat \Sigma^{(r+1)}=\hat
  \Sigma^{(r,p)}$. Increment the 
counter $r$ to $r+1$. Go to 2. 
\end{enumerate}
\end{algorithm}
The iterations can be stopped according to a 
criterion such as ``the estimate of
$\Sigma$ is not changed'' (in some pre-determined accuracy).

\begin{example}
\label{example}
\rm
Figure \ref{fig:pseudo-var} illustrates ICF for the 
model $\mathbf{N}(G)$ based on the
graph $G$ shown in Figure \ref{fig:ex.bi}. 
The algorithm cycles in arbitrary order through the
four regressions 
$(Y_i\mid Y_{-i})$, $i=1,2,3,4$. 
In Figure \ref{fig:pseudo-var}, a filled circle represents
variables in the conditioning set $-i$, and an unfilled circle stands for
the variable $i$ forming the response variable in the considered
regression. The directed edges coincide with bi-directed edges in
the original graph in Figure \ref{fig:ex.bi} and indicate the
pseudo-variable regressions to be carried out.  The vertices that are
joined to vertex $i$ by a directed edges are labelled with the 
pseudo-variables that act as covariates. The directed edges are
labelled with the covariances that are estimated.  
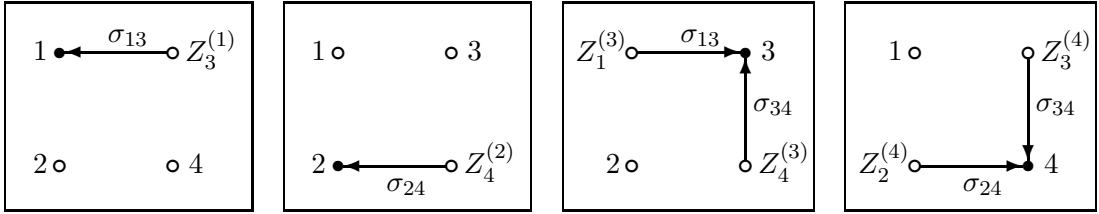
\begin{figure}[htb]
\bigskip
\begin{center}
\setlength{\unitlength}{0.29pt}
\fbox{
\begin{picture}(275,250)(10,-10)
\thicklines
\linethickness{0.65pt}
\put(32,176){\makebox(0,0)[b]{1}}
\put(32,28){\makebox(0,0)[b]{2}}
\put(255,164){\makebox(0,0)[b]{$Z^{(1)}_3$}}
\put(235,28){\makebox(0,0)[b]{4}}
\put(197.9,185.6){\vector(-1,0){   134.0}}
\put(120,200){$\sigma_{13}$}
\put(56.9,185.6){\circle*{14}}
\put(56.9, 37.6){\circle{14}}
\put(204.9,185.6){\circle{14}}
\put(204.9, 37.6){\circle{14}}
\end{picture} }
\hspace{0.15cm}
\fbox{
\begin{picture}(275,250)(10,-10)
\thicklines
\linethickness{0.65pt}
\put(32,176){\makebox(0,0)[b]{1}}
\put(32,28){\makebox(0,0)[b]{2}}
\put(235,176){\makebox(0,0)[b]{3}}
\put(255,17){\makebox(0,0)[b]{$Z^{(2)}_4$}}
\put(120,5){$\sigma_{24}$}
\put(197.9, 37.6){\vector(-1,0){   134.0}}
\put(56.9,185.6){\circle{14}}
\put(56.9, 37.6){\circle*{14}}
\put(204.9,185.6){\circle{14}}
\put(204.9, 37.6){\circle{14}}
\end{picture} }
\hspace{0.15cm}
\fbox{
\begin{picture}(280,250)(-10,-10)
\thicklines
\linethickness{0.65pt}
\put(17,164){\makebox(0,0)[b]{$Z^{(3)}_1$}}
\put(32,28){\makebox(0,0)[b]{2}}
\put(235,176){\makebox(0,0)[b]{3}}
\put(255,17){\makebox(0,0)[b]{$Z^{(3)}_4$}}
\put(63.9,185.6){\vector(1,0){   134.0}}
\put(204.9,44.6){\vector(0,1){   134.0}}
\put(120,200){$\sigma_{13}$}
\put(215,110){$\sigma_{34}$}
\put(56.9,185.6){\circle{14}}
\put(56.9, 37.6){\circle{14}}
\put(204.9,185.6){\circle*{14}}
\put(204.9, 37.6){\circle{14}}
\end{picture} }
\hspace{0.15cm}
\fbox{
\begin{picture}(280,250)(-10,-10)
\thicklines
\linethickness{0.65pt}
\put(32,176){\makebox(0,0)[b]{1}}
\put(17,17){\makebox(0,0)[b]{$Z^{(4)}_2$}}
\put(255,164){\makebox(0,0)[b]{$Z^{(4)}_3$}}
\put(235,28){\makebox(0,0)[b]{4}}
\put(120,5){$\sigma_{24}$}
\put(215,110){$\sigma_{34}$}
\put(63.9, 37.6){\vector(1,0){   134.0}}
\put(204.9,178.6){\vector(0,-1){   134.0}}
\put(56.9,185.6){\circle{14}}
\put(56.9, 37.6){\circle{14}}
\put(204.9,185.6){\circle{14}}
\put(204.9, 37.6){\circle*{14}}
\end{picture} }
  \caption{Illustration of the pseudo-variable regressions in ICF.}
  \label{fig:pseudo-var}
\end{center}
\end{figure}  
\end{example}

\begin{remark}[Complexity]\rm 
The algorithm  can be restated only in terms of the empirical
covariance matrix $S$ defined in \eqref{empcov}. 
For example in \eqref{ls},
\begin{equation}
  \label{eq:expressS}
  \begin{split}
    Y_i (Z^{(i)}_{\spo(i)})'&=
    S_{i,-i}[(\Sigma_{-i,-i})^{-1}]_{-i,\spo(i)},\\
    Z^{(i)}_{\spo(i)}(Z^{(i)}_{\spo(i)})'&=
    [(\Sigma_{-i,-i})^{-1}]_{\spo(i),-i}S_{-i,-i}
    [(\Sigma_{-i,-i})^{-1}]_{-i,\spo(i)}. 
  \end{split}
\end{equation}
Other products between data matrices appearing in the sequel can be
similarly expressed in terms of the empirical covariance matrix $S$. 
Thus, the sample size does
not affect the complexity of the algorithm.
The complexity of one of the algorithm's pseudo-variable regression
steps is dominated by the computation of the inverse of
$\Sigma_{-i,-i}$  in (\ref{eq:2}), and the inversion of a matrix of
size $\spo(i)\times \spo(i)$ \eqref{ls}. Note that $\Sigma_{-i,-i}$
may be sparse and special methods for inversion of sparse matrices
might be useful. In particular, if the induced subgraph $G_{-i}$ has
disconnected components then only the submatrices of $\Sigma$ over
connected components containing spouses of $i$ have to be inverted.
\end{remark}

\subsection{Convergence}\label{convergence} 

The key to prove convergence of ICF is to recognize that the algorithm
consists of iterated partial maximizations over sections of the
parameter space $\mathbf{P}(G)$.  In ICF we repeatedly maximize
the likelihood function of the covariance graph model partially by
allowing only the entries in the $i$-th row and column of $\Sigma$ to
vary. The remaining entries are fixed.  A bit more formally, we
consider the parameter space  
\begin{equation}
\Theta= \{\Sigma\in\mathbf{P}(G)\mid \ell(\Sigma)\ge
\ell(\hat\Sigma^{(0)})\},
\end{equation}
which is compact, though not necessarily connected, and contains the
global maximizer of $\ell(\Sigma)$. Recall that we assume the
empirical covariance matrix $S$ to be positive definite.
Defining the section $\Theta_i(\bar\Sigma)\subsetneq \Theta$ as
\begin{equation}
\Theta_i(\bar\Sigma) = \big\{\Sigma\in\Theta \mid  \Sigma_{-i,-i}=
\bar\Sigma_{-i,-i}\big\},
\end{equation}
it becomes clear that the algorithm steps
2(i)-2(iii) maximize the log-likelihood function partially over the section
$\Theta_i(\hat\Sigma^{(r,i-1)})$, i.e. 
\begin{equation}\label{maxsect}
\hat\Sigma^{(r,i)}=\arg\max \big\{ \ell(\Sigma)\mid
\Sigma\in\Theta_i(\hat\Sigma^{(r,i-1)})\big\}. 
\end{equation}
This local and global maximizer over the section is unique.  If a
matrix $\Sigma\in\mathbf{P}(G)$ maximizes the log-likelihood function
over all sections $\Theta_i(\Sigma)$, $i \in V$, simultaneously, then
it solves the likelihood equations.  Hence, the following theorem
follows from results in \citet[Appendix]{drtoneichler:2004}. 
\begin{theorem}
\label{conv}
Suppose the sequence $(\hat\Sigma^{(r)})$ is constructed by 
the ICF algorithm. Then all accumulation points
of $(\hat\Sigma^{(r)})$ are saddle points or local maxima 
of the log-likelihood function. Moreover, all accumulation points have the same
likelihood value.  In particular, if the likelihood equations have only
finitely many solutions, then $(\hat\Sigma^{(r)})$ converges.  
\end{theorem}

\section{Iterative conditional fitting with multivariate updates}
\label{sec:icfsur}

The algorithm presented in \S \ref{newmethod} is based on updating one
row and column of an estimate of the covariance matrix
$\Sigma\in\mathbf{P}(G)$ by carrying out a univariate regression.  A
natural modification of this 
approach is to update several rows and columns of the estimate 
$\Sigma\in\mathbf{P}(G)$ simultaneously using multivariate
regression. 

\subsection{Seemingly unrelated pseudo-variable regressions}
\label{sec:pseudosur}

Let $C\subseteq V$ be a subset of the vertices. In order to estimate
all rows and columns of $\Sigma$ that are indexed by the vertices in $C$
in the ICF algorithm presented in \S \ref{newmethod}, we have to
carry out several univariate pseudo-variable regressions for $(Y_i\mid
Y_{-i})$, $i\in C$. Instead, we would like to 
consider only one multivariate regression of the form $(Y_C\mid
Y_{-C})$, where $-C=V\setminus C$.  
The conditional distribution  
\begin{equation}
\label{eq:conddistC}
(Y_C\mid Y_{-C}) \sim \ND_{C\times
  N}(B_CY_{-C}, \Lambda_C\otimes I_N)
\end{equation}
is specified by 
the matrix of
regression coefficients
\begin{equation}
  \label{eq:BC}
  B_C=\Sigma_{C,-C}(\Sigma_{-C,-C})^{-1}\in \RRR^{C\times -C},
\end{equation}
and the conditional covariance matrix
\begin{equation}
  \label{eq:lambdaC}
  \Lambda_C=\Sigma_{C,C}-\Sigma_{C,-C}(\Sigma_{-C,-C})^{-1}\Sigma_{-C,C}\in
  \mathbf{P}(C).
\end{equation}
In order for the conditional distribution (\ref{eq:conddistC}) to be
of a simple structure, there
should be no constraints on the $\Lambda_{C}$, in which case
$\mathbf{P}(C)=\mathbf{P}_C(G)$. This holds if 
there are no constraints on the submatrix $\Sigma_{C,C}$, which in
turn holds if the set $C$ is complete, i.e., if  $i\bi j$
whenever $i,j\in C$ and $i\not= j$.  Then the only constraints on the
conditional distribution  (\ref{eq:conddistC}) are on the
matrix of regression coefficients $B_C$ and stem from restrictions
that $\sigma_{ij}=0$, if $i\in C$, 
$j\not\in C$ and $j\not\bi i$.

Let 
\begin{equation}
  \spo(C)=\big[\cup (\spo(i)\mid i\in C)\big]\setminus C
\end{equation} 
be the {\it spouses of}
$C$, that is the vertices that are not in $C$ but adjacent to some
vertex in $C$, and let $\nsp(C)=V\setminus( \spo(C)\cup C)$ be the
{\it non-spouses of} $C$, yielding 
the partition $V=C\cup \spo(C)\cup \nsp(C)$.
If we define the pseudo-variables
\begin{equation}
\label{eq:ZC}
 Z^{(C)}_{\spo(C)} = [(\Sigma_{-C,-C})^{-1}]_{\spo(C),-C}\,Y_{-C}
 \in\RRR^{\spo(C)\times N},
\end{equation}
then we can rewrite (\ref{eq:conddistC}) as
\begin{equation}
\label{eq:pseudoC}
(Y_C\mid Y_{-C}) \sim \ND_{C\times
  N}(\Sigma_{C,\spo(C)} \,Z^{(C)}_{\spo(C)}, \Lambda_C\otimes I_N),
\end{equation}
because $\Sigma_{C,\nsp(C)}=0$.
As $\Sigma$ ranges through $\mathbf{P}(G)$, the submatrix
$\Sigma_{C,\spo(C)}$ playing the role of regression coefficients in
(\ref{eq:pseudoC}) ranges through the linear space
\begin{equation}
  \label{eq:sigmaCrange}
  \mathbf{P}_{C,\spo(C)}(G) = \big\{ A\in \RRR^{C\times \spo(C)}\mid A_{ij}=0
  \mbox{ if } i\not\bi j\big\}.
\end{equation}
Hence, (\ref{eq:pseudoC}) constitutes seemingly unrelated
regressions \citep{zellner:62}.

\subsection{The iterative conditional fitting algorithm with multivariate updates}
\label{sec:estpseudosur}

ML estimation in seemingly unrelated regressions itself generally
requires iterative algorithms, such as iterating the two-step
estimator of \cite{zellner:62}.  In the case of (\ref{eq:pseudoC}),
the two-step estimator consists of first estimating $\Sigma_{C,\spo(C)}$
for some fixed $\Lambda_C$ by generalized least squares, and then
estimating $\Lambda_C$ as the empirical covariance matrix of the
residuals $Y_i-\Sigma_{C,\spo(C)}Z^{(C)}_{\spo(C)}$ computed with the
estimate of $\Sigma_{C,\spo(C)}$ obtained in
the first step. However, if the 
current estimate of $\Sigma$ is used to obtain starting values
$\Sigma_{C,\spo(C)}$ and $\Lambda_C$, then the two-step method does not
have to be iterated in order to obtain estimates for the seemingly unrelated
pseudo-regressions (\ref{eq:pseudoC}) that yield a convergent
ICF algorithm with multivariate updates. 
For specification of the estimator of $\Sigma_{C,\spo(C)}$ we need to
introduce the matrix $P_C$ of the linear map that sends the
vector of unrestricted elements in 
$\Sigma_{C,\spo(C)}$ to the matrix
$\Sigma_{C,\spo(C)}\in\mathbf{P}_{C,\spo(C)}(G)$.  The
vector of unrestricted elements of $\Sigma_{C,\spo(C)}$ is the vector
$\sigma_{C}=(\sigma_{ij}\mid i\in C,\; j\in \spo(C),\; i\bi j)$.  
The matrix $P_C$ has 
exactly one entry equal to one in each column, the other entries are
zero, and it satisfies $\cvec(\Sigma_{C,\spo(C)})=P_C \sigma_C$ for
$\Sigma\in \mathbf{P}(G)$;
compare the definition of the matrix $Q$ in \S\ref{covgraphmodel}.

In order to run ICF with multivariate updates, we have to choose a
family of complete sets $(C\mid C\in \mathcal{C})$ such that 
\begin{equation}
  \label{eq:chooseC}
  \cup (C \mid C\in\mathcal{C}) =  V,
\end{equation}
where the sets $C$ do not have to be disjoint.  For example the sets
$C$ could be chosen as edges, but the largest possible choice for the
sets $C$ would be the cliques, i.e., the maximal complete sets, in
$G$.

\begin{algorithm}
  \label{alg:icfsur}
  For a given choice of $\mathcal{C}$, the ICF
  algorithm with multivariate updates can be implemented as:
\begin{enumerate}
\item[1.] {\em (Initialization)} Set the iteration counter $r=0$, and
  choose a starting value $\hat 
\Sigma^{(0)}\in\mathbf{P}(G)$, e.g.\ the identity matrix $\hat
\Sigma^{(0)}=I_V$.  
\item[2.] {\em (Updates)} Order the sets in the family $\mathcal{C}$ as 
$\mathcal{C}=\{C_1,\ldots,C_q\}$, set $\hat \Sigma^{(r,0)}=\hat \Sigma^{(r)}$, 
and repeat the following steps for all $C_k\in \mathcal{C}$:
\begin{enumerate}
\item[(i)] Let $\hat \Sigma^{(r,k)}_{-{C_k},-{C_k}}=\hat
  \Sigma^{(r,k-1)}_{-{C_k},-{C_k}}$. From this submatrix, compute
  the conditional covariance matrix 
  $\hat\Lambda_{C_k}$ according to (\ref{eq:lambdaC})  and
  the pseudo-variables $Z^{(k)}_{\spo(C_k)}$ according to
  (\ref{eq:ZC}). Calculate $\hat\Omega_{C_k}=(\hat\Lambda_{C_k})^{-1}$.
\item[(ii)] Compute the (generalized least
  squares) matrix that satisfies $\cvec (\hat
  \Sigma^{(r,k)}_{C_k,\spo(C_k)}) = P_{C_k}\hat\sigma_{C_k}$, where 
\begin{multline}\label{lsC1}
  \hat\sigma_{C_k} =
  \Big\{ P_{C_k}'\,\big\{\,
  [Z^{(k)}_{\spo(C_k)}(Z^{(k)}_{\spo(C_k)})']\otimes  
  \hat\Omega_{C_k}\big\}\,P_{C_k} \Big\}^{-1} \times \\ \big\{
  P_{C_k}'\,\cvec[\hat\Omega_{C_k}\,
  Y_C\,(Z^{(k)}_{\spo(C_k)})']\big\}. 
\end{multline}
\item[(iii)] Compute the empirical covariance matrix
  of residuals
\begin{equation}
\label{eq:lsC2}
\hat\Lambda_{C_k}=
\frac{1}{n}(Y_{C_k}-\hat
  \Sigma^{(r,k)}_{C_k,\spo(C_k)}Z^{(k)}_{\spo(C_k)})(Y_{C_k}-\hat
  \Sigma^{(r,k)}_{C_k,\spo(C_k)}Z^{(k)}_{\spo(C_k)})'.
\end{equation}
\item[(iii)] Complete $\hat\Sigma^{(r,k)}$ by
setting 
\begin{equation}
  \hat\Sigma^{(r,k)}_{C_k,C_k}=\hat\Lambda_{C_k}+\hat
  \Sigma^{(r,k)}_{C_k,\spo(C_k)} \big[(\hat
    \Sigma^{(r,k)}_{-{C_k},-{C_k}})^{-1}\big]_{\spo(C_k),\spo(C_k)}
      \hat\Sigma^{(r,k)}_{\spo(C_k),C_k};
\end{equation}
compare (\ref{eq:lambdaC}). 
\end{enumerate}
\item[3.] {\em (Repeat)} Set $\hat \Sigma^{(r+1)}=\hat
  \Sigma^{(r,q)}$. Increment the 
counter $r$ to $r+1$. Go to 2. 
\end{enumerate}
\end{algorithm}
Note that if the family $\mathcal{C}$ consists of only singletons then
Algorithm \ref{alg:icfsur} reduces to Algorithm \ref{alg:icf}.

\begin{example}
\label{sec:exicfsur}
  \rm
We take up the covariance graph shown in Figure \ref{fig:ex.bi}.  For the
family $\mathcal{C}$ of complete vertex sets, several choices are possible.
If the cliques $\mathcal{C}=\{13,34,24\}$ are chosen, then all conditional
distributions considered in ICF are bivariate, whereas for
$\mathcal{C}=\{1,2,34\}$ two univariate distributions are estimated in
conjunction with a bivariate distribution.  For the clique choice
$\mathcal{C}=\{13,34,24\}$, we illustrate the seemingly unrelated
pseudo-variable regressions to be estimated in Figure \ref{fig:ex.icfsur},
which is to be interpreted similarly as Figure \ref{fig:pseudo-var}.  An
additional feature are the bi-directed edges that connect the vertices in
the sets $C\in\mathcal{C}$; see \citet{richardson:2002} for a formal
definition of these graphs.

\begin{figure}[htb]
\bigskip
\begin{center}
\setlength{\unitlength}{0.29pt}
\fbox{
\begin{picture}(290,250)(0,-10)
\thicklines
\linethickness{0.65pt}
\put(32,176){\makebox(0,0)[b]{1}}
\put(32,28){\makebox(0,0)[b]{2}}
\put(235,176){\makebox(0,0)[b]{3}}
\put(265,17){\makebox(0,0)[b]{$Z^{(13)}_4$}}
\put(197.9,185.6){\vector(-1,0){   134.0}}
\put(63.9,185.6){\vector(1,0){   134.0}}
\put(204.9,44.6){\vector(0,1){   134.0}}
\put(215,110){$\sigma_{34}$}
\put(120,200){$\sigma_{13}$}
\put(56.9,185.6){\circle*{14}}
\put(56.9, 37.6){\circle{14}}
\put(204.9,185.6){\circle*{14}}
\put(204.9, 37.6){\circle{14}}
\end{picture} }
\hspace{0.25cm}
\fbox{
\begin{picture}(265,250)(-20,-10)
\thicklines
\linethickness{0.65pt}
\put(12,165){\makebox(0,0)[b]{$Z^{(34)}_1$}}
\put(12,17){\makebox(0,0)[b]{$Z^{(34)}_2$}}
\put(235,176){\makebox(0,0)[b]{3}}
\put(235,28){\makebox(0,0)[b]{4}}
\put(63.9,185.6){\vector(1,0){   134.0}}
\put(204.9,44.6){\vector(0,1){   134.0}}
\put(63.9, 37.6){\vector(1,0){   134.0}}
\put(204.9,178.6){\vector(0,-1){   134.0}}
\put(120,200){$\sigma_{13}$}
\put(120,5){$\sigma_{24}$}
\put(215,110){$\sigma_{34}$}
\put(56.9,185.6){\circle{14}}
\put(56.9, 37.6){\circle{14}}
\put(204.9,185.6){\circle*{14}}
\put(204.9, 37.6){\circle*{14}}
\end{picture} }
\hspace{0.25cm}
\fbox{
\begin{picture}(290,250)(10,-10)
\thicklines
\linethickness{0.65pt}
\put(32,176){\makebox(0,0)[b]{1}}
\put(32,28){\makebox(0,0)[b]{2}}
\put(265,165){\makebox(0,0)[b]{$Z^{(24)}_3$}}
\put(235,28){\makebox(0,0)[b]{4}}
\put(120,5){$\sigma_{24}$}
\put(215,110){$\sigma_{34}$}
\put(63.9, 37.6){\vector(1,0){   134.0}}
\put(204.9,178.6){\vector(0,-1){   134.0}}
\put(197.9, 37.6){\vector(-1,0){   134.0}}
\put(56.9,185.6){\circle{14}}
\put(56.9, 37.6){\circle*{14}}
\put(204.9,185.6){\circle{14}}
\put(204.9, 37.6){\circle*{14}}
\end{picture} }
  \caption{Illustration of the seemingly unrelated pseudo-variable
    regressions in ICF with multivariate updates, and $\mathcal{C}=\{13,34,24\}$.} 
  \label{fig:ex.icfsur}
\end{center}
\end{figure}
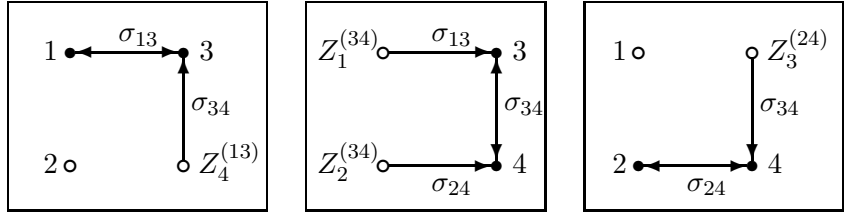
\end{example}

\subsection{Convergence}

The ICF algorithm with multivariate updates is still an iterative
partial maximization algorithm.  However, the sections in the
parameter space over which maximizations are performed are not quite
as simple as the sections described in \S \ref{convergence}.
Steps 2(ii) and 2(iii) of Algorithm \ref{alg:icfsur} do not jointly
maximize the log-likelihood function $\ell$ over sections of
the form 
\begin{equation}
  \label{eq:sectC}
  \Theta_C(\bar\Sigma) = \big\{\Sigma\in\Theta \mid  \Sigma_{-C,-C}=
  \bar\Sigma_{-C,-C}\big\}.
\end{equation}
Instead step 2(ii) maximizes $\ell$ over sections of the form
\begin{equation}
  \label{eq:sectC1}
  \Theta_{1,C}(\bar\Sigma) = \big\{\Sigma\in\Theta \mid  \Sigma_{-C,-C}=
  \bar\Sigma_{-C,-C},\; \Lambda_C= \bar\Lambda_C\big\},
\end{equation}
where $\Lambda_C$ is again the conditional covariance matrix from
(\ref{eq:lambdaC}).  The subsequent step 2(iii) maximizes $\ell$ over
sections of the form  
\begin{equation}
  \label{eq:sectC2}
  \Theta_{2,C}(\bar\Sigma) = \big\{\Sigma\in\Theta \mid  \Sigma_{-C,-C}=
  \bar\Sigma_{-C,-C},\; \Sigma_{C,-C}= \bar\Sigma_{C,-C}\big\}.
\end{equation}
Nevertheless it holds under condition (\ref{eq:chooseC}) that if
$\Sigma$ maximizes the log-likelihood function $\ell$ over both  
section $\Theta_{1,C}(\bar\Sigma)$ and $\Theta_{2,C}(\bar\Sigma)$
simultaneously for all $C\in \mathcal{C}$, then $\Sigma$ is a solution
to the likelihood equations.  Thus, Theorem \ref{conv} holds also for
ICF with multivariate updates as stated in Algorithm \ref{alg:icfsur}.

\section{Empirical likelihood estimation}
\label{sec:empir}

In contexts where it is not appropriate to assume multivariate
normality, we may still wish to estimate a covariance matrix subject
to various zero restrictions. Here we present an approach based on
empirical likelihood \citep{owen:2001}.  In the resulting method an
estimate of the underlying distribution is obtained by maximizing a
non-parametric likelihood under constraints that include the desired
zero covariance restrictions; see \citet*{scmsh1} and
\citet*{hellimbens} for similar applications of empirical likelihood.

We associate a weight $w_k$ with the $k$-th sample observation
$Y^{(k)}_V$, $k\in N$. Estimating the mean vector and covariance
matrix simultaneously, we solve the nested constrained maximization
problem
\begin{equation}
\label{eq:objecfun}
\max_{{\mu}}\left\{\raisebox{-.6ex}{$\substack{\hbox{max}\\{w}
=\left(w_1,\ldots,w_n\right)}$}\prod_{k\in N}n w_k\right\}
\end{equation}
subject to
\begin{align}
  ~&w_k\ge 0, k\in N,\label{eq:nneg}\\
  &\sum_{k\in N}w_k=1,\label{eq:sumto1}\\
  &\sum_{k\in N}w_k \left (Y^{(k)}_i-\mu_i\right)=0,
  \quad\forall~i\in V,\label{eq:mean0}\\
  &\sum_{k\in N}w_k
  \left(Y^{(k)}_i-\mu_i\right)\left(Y^{(k)}_j-\mu_j\right)=0,
  \quad\text{$\forall~ i,j\in V$ s.t. $i\not\leftrightarrow
    j$}\label{eq:cov0}.
\end{align}
Without the additional constraints (\ref{eq:mean0}) and
(\ref{eq:cov0}), the empirical likelihood ratio $\prod_{k\in N} n w_k$ is
maximized for $w_k=1/n$, $k\in N$.  The additional constraint
\eqref{eq:mean0} enforces that the mean of the reweighted rows of $Y$
is equal to ${\mu}$.  Constraint \eqref{eq:cov0} ensures that the
estimated weights $\h{w}_k$ are such that the empirical covariance
matrix of the reweighted sample satisfies the zero constraints
specified by the graph $G$.

In order to avoid obvious problems with feasibility of the
optimization problem, we assume that the sample size, i.e., the number
of weights, is strictly larger than the number of constraints in
\eqref{eq:sumto1} and \eqref{eq:cov0}.  Note that the number of
constraints \eqref{eq:cov0} may grow quadratically as the number of
variables increases.  The nesting of the maximization steps in
\eqref{eq:objecfun} is done to avoid cubic 
constraints in $w$, which would have resulted had we substituted
\begin{equation}
\mu_i=\sum_{k\in N}w_k Y^{(k)}_i,\quad \forall i\in V
\end{equation}
in \eqref{eq:cov0} and made the constraints in \eqref{eq:mean0}
redundant.  The constrained maximization problem can be solved through
its dual problem, in which the number of unknowns is equal to the
number of constraints of the original problem; see \citet{owen:2001}
and \citet{scmsh1} for details.
 
If ${\h{\mu}}$ and ${\h{w}}$ are, respectively, the vectors of mean and
weights maximizing \eqref{eq:objecfun} under the constraints
\eqref{eq:nneg}-\eqref{eq:cov0}, then the estimated covariance matrix is
given by
\begin{equation}
\h{\Sigma}_{E}=\left(Y-{\h{\mu}}\otimes 1_{N}\right)\cdot 
\hbox{diag}({\h{w}})\cdot \left(Y-{\h{\mu}}\otimes 1_{N}\right)', 
\end{equation}
where $\hbox{diag}({\h{w}})$ is an $n\times n$ diagonal matrix with
${\h{w}}$ along its diagonal.  Following \citet{owen:2001} and
\citet{Qinlawless1} one can show that asymptotically ${\h{\mu}}$ and
$\h{\Sigma}_E$ are consistent.  


\section{Data and simulations} 
\label{ch4data}

We now compare the three approaches to estimation of a covariance matrix
with zeros in a data example and in simulations: (i) ML estimation relying
on ICF, (ii) dual likelihood estimation as described in \S \ref{sec:dual},
and (iii) empirical likelihood estimation.

\subsection{Gene expression in yeast}
\label{sec:yeast}

\citet{gasch2000} present gene expression data from microarray
experiments with yeast strands.  We focus on $p=8$ genes related to
galactose utilization.  The gene {\em GAL11} is responsable for
transcription.  The genes {\em GAL4} and {\em GAL80} are involved in
galactose regulation.  Gene {\em GAL2} is related to transport and the
remaining four genes, {\em GAL1}, {\em GAL3}, {\em GAL7}, and {\em
  GAL10}, are involved in galactose metabolism.  There are $n=134$
experiments with gene expression measurements for all eight genes.
The observed marginal correlations and standard deviations are shown
in Table \ref{table1}, where we denote the variables for the gene
expression measurements by $X_i$, $i=1,2,3,4,7,10,11,80$, using the
obvious correspondence.

\begin{table}[htbp]
\caption{\label{table1}Observed marginal correlations and standard
deviations.}  
\medskip
 \centering\small
\begin{tabular}{lrrrrrrrr}
& $X_{11}$ &  $X_{4}$ & $X_{80}$ & $X_{2}$ & $X_{1}$ & $X_{3}$ & $X_{7}$ & $X_{10}$\\
$X_{11}$ &  \\ 
$X_{4}$  &  0.24 &  \\  
$X_{80}$ &  0.08 &  0.23 &  \\  
$X_{2}$  &  $-$0.18 & $-$0.03 &  0.26 &  \\  
$X_{1}$  &  $-$0.10 & $-$0.10 &  0.28 &  0.87 &  \\
$X_{3}$  &  $-$0.18 &  0.12 &  0.20 &  0.44 &  0.39 &  \\
$X_{7}$  &  $-$0.07 & $-$0.08 &  0.21 &  0.81 &  0.88 &  0.50 &  \\
$X_{10}$ & $-$0.08 & $-$0.07 &  0.26 &  0.87 &  0.92 &  0.46 &  0.91 &  \\[0.1cm]
SD &  0.39 & 0.36 &  0.47 & 1.70 & 1.70 & 0.78 & 1.85 & 1.54 
\end{tabular}
\end{table}

\begin{figure}[t]
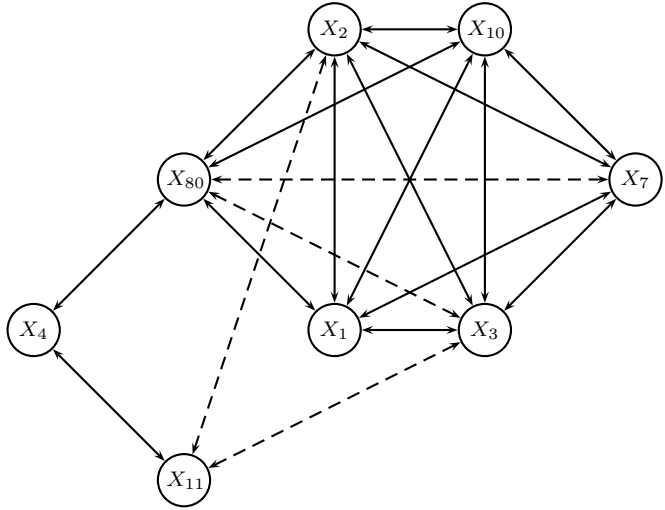

    \vspace{9cm}\hspace{1cm}
    \scriptsize
    \psset{linewidth=0.8pt}          
    \settowidth{\MyLength}{$X_{0}$}
    \newcommand{\myNode}[2]{\circlenode{#1}{\makebox[\MyLength]{#2}}}
    \rput(2, 4){\myNode{1}{$X_{4}$}} 
    \rput(4, 2){\myNode{2}{$\,X_{11}$}}
    \rput(4, 6){\myNode{3}{$\,X_{80}$}}
    \rput(6, 4){\myNode{4}{$X_1$}} 
    \rput(6, 8){\myNode{5}{$X_2$}}
    \rput(8, 4){\myNode{6}{$X_3$}}
    \rput(8, 8){\myNode{7}{$\,X_{10}$}}
    \rput(10, 6){\myNode{8}{$X_7$}} 
    \ncline{<->}{1}{2}
    \ncline{<->}{1}{3}
    \ncline{<->}{3}{4}
    \ncline{<->}{3}{5}
    \ncline{<->}{3}{7}
    \ncline{<->}{4}{5}
    \ncline{<->}{4}{6}
    \ncline{<->}{4}{7}
    \ncline{<->}{4}{8}
    \ncline{<->}{5}{6}
    \ncline{<->}{5}{7}
    \ncline{<->}{5}{8}
    \ncline{<->}{6}{7}
    \ncline{<->}{6}{8}
    \ncline{<->}{7}{8}
    \psset{linestyle=dashed}          
    \ncline{<->}{2}{5}
    \ncline{<->}{2}{6}
    \ncline{<->}{3}{6}
    \ncline{<->}{3}{8}
    \vspace{-1.5cm}
\caption{\label{fig:SIN} Covariance graph for data in Table
  \ref{table1}.}
\end{figure}

By multiple testing of correlations as described in
\cite{drtonperlman:2004,drtonperlman:2004b} and implemented in the R
package `SIN', we selected the two covariance graphs $G_\mathrm{s}\subset
G_\mathrm{d}$ that are illustrated in Figure \ref{fig:SIN}.  The larger
graph $G_\mathrm{d}$ contains all edges shown, i.e., both the solid and the
dashed edges, whereas the sub-graph $G_\mathrm{s}$ includes only the solid
edges.  In $G_\mathrm{s}$ the vertices 1, 2, 3, 7, and 10 form a clique and
in $G_\mathrm{d}$ the clique is enlarged to include vertex 80.  With the R
package `ggm' and additional code, we computed the ML, the dual, and the
empirical likelihood estimates of the covariance matrix under the zero
constraints specified in the graphs.  The results for both $G_\mathrm{s}$
and $G_\mathrm{d}$ are shown in Table \ref{tableGs}. We remark that we
started ICF from the identity matrix and that Anderson's algorithm gave the
same results as ICF.  However, although we refer to ``ML estimates'', ICF
is only guaranteed to find a stationary point which may not be the global
maximizer of the likelihood.

\begin{table}[htbp]
\caption{\label{tableGs} Marginal correlations and standard deviations
from ML (M), dual (D), and empirical likelihood (E) estimates for graph
$G_\mathrm{s}$ (lower half) and graph $G_\mathrm{d}$ (upper italicized
half).}  
\medskip
 \centering\small
\begin{tabular}{lrrrrrrrrrc}
&  $X_{11}$ &  $X_{4}$ & $X_{80}$ & $X_{2}$ & $X_{1}$ & $X_{3}$ &
$X_{7}$  & $X_{10}$ & SD\\
$X_{11}$& &\it 0.28& \it 0& \it $-$0.12& \it 0& \it $-$0.21& \it 0& \it 0 & \it 0.40&M\\ 
        & &\it 0.26& \it 0& \it $-$0.11& \it 0& \it $-$0.20& \it 0& \it 0 & \it 0.39&D\\ 
        & &\it 0.25& \it 0& \it $-$0.11& \it 0& \it $-$0.20& \it 0& \it 0& \it 0.39&E\\[0.1cm] 
$X_{4}$ & 0.22 &  & \it 0.20& \it 0& \it 0& \it 0& \it 0& \it 0& \it 0.36&M\\  
        & 0.27 &  & \it 0.21& \it 0& \it 0& \it 0& \it 0& \it 0& \it 0.35&D\\  
        & 0.28 &  & \it 0.27& \it 0& \it 0& \it 0& \it 0& \it 0& \it 0.36&E\\[0.1cm] 
$X_{80}$& 0 & 0.22 &   & \it 0.27& \it 0.29& \it 0.19& \it 0.22& \it 0.27& \it 0.47&M\\  
        & 0 & 0.20 &   & \it 0.28& \it 0.31& \it 0.19& \it 0.23& \it 0.28& \it 0.47&D\\  
        & 0 & 0.18 &   & \it 0.26& \it 0.31& \it 0.16& \it 0.21& \it 0.27& \it 0.48&E\\[0.1cm]   
$X_{2}$ & 0 & 0 &  0.08 &  & \it 0.86& \it 0.43& \it 0.81& \it 0.87& \it 1.69&M\\  
        & 0 & 0 &  0.09 &  & \it 0.86& \it 0.43& \it 0.81& \it 0.87& \it 1.68&D\\  
        & 0 & 0 &  0.17 &  & \it 0.83& \it 0.43& \it 0.79& \it 0.85& \it 1.48&E\\[0.1cm]  
$X_{1}$ & 0 & 0 &  0.11 &  0.86 &   & \it 0.38& \it 0.88& \it 0.92& \it 1.70&M\\
        & 0 & 0 &  0.12 &  0.86 &   & \it 0.39& \it 0.88& \it 0.91& \it 1.69&D\\
        & 0 & 0 &  0.10 &  0.83 &   & \it 0.34& \it 0.85& \it 0.88& \it 1.48&E\\[0.1cm] 
$X_{3}$ & 0 & 0 &  0 &  0.43 &  0.38 &   & \it 0.49& \it 0.44& \it 0.78&M\\
        & 0 & 0 &  0 &  0.39 &  0.37 &   & \it 0.51& \it 0.46& \it 0.78&D\\
        & 0 & 0 &  0 &  0.39 &  0.31 &   & \it 0.49& \it 0.46& \it 0.78&E\\[0.1cm] 
$X_{7}$ & 0 & 0 &  0 &  0.81 &  0.88 &  0.50 &   & \it 0.91& \it 1.85&M\\
        & 0 & 0 &  0 &  0.80 &  0.87 &  0.50 &   & \it 0.91& \it 1.84&D\\
        & 0 & 0 &  0 &  0.77 &  0.83 &  0.38 &   & \it 0.90& \it 1.68&E\\[0.1cm] 
$X_{10}$& 0 & 0 &  0.08 &  0.87 &  0.91 &  0.45 &  0.91 &  & \it 1.54&M\\
        & 0 & 0 &  0.08 &  0.86 &  0.91 &  0.44 &  0.90 &  & \it 1.53&D\\
        & 0 & 0 &  0.13 &  0.86 &  0.87 &  0.36 &  0.88 &  & \it 1.36&E\\[0.3cm]
SD      & 0.39 & 0.36 &  0.47 & 1.70 & 1.70 & 0.78 & 1.85 & 1.54&&M\\ 
        & 0.37 & 0.35 &  0.45 & 1.61 & 1.61 & 0.75 & 1.79 & 1.47&&D\\
        & 0.38 & 0.33 &  0.47 & 1.41 & 1.37 & 0.74 & 1.57 & 1.22&&E
\end{tabular}
\end{table}

Upon inspection of Table \ref{tableGs} we find that the three estimates are
in better agreement for the graph $G_\mathrm{d}$.  This graph yields the
better fitting covariance graph model.  The deviance of the model
$\mathbf{N}(G_\mathrm{d})$ under comparison to the model based on the
complete graph equals 9.98 over 9 degrees of freedom, whereas the deviance
of the model $\mathbf{N}(G_\mathrm{s})$ equals 33.07 over 13 degrees of
freedom.  This indicates a good fit of $\mathbf{N}(G_\mathrm{d})$ and a
poor fit of the more restrictive model $\mathbf{N}(G_\mathrm{s})$.  The
difference in log-likelihood between ML and dual estimates equals 4.29 in
$\mathbf{N}(G_\mathrm{s})$ and reduces to 0.51 in
$\mathbf{N}(G_\mathrm{d})$.  In contrast the difference in log-likelihood
between ML and empirical likelihood estimates equals 20.54 in
$\mathbf{N}(G_\mathrm{s})$ and 5.67 in $\mathbf{N}(G_\mathrm{d})$.


\subsection{Simulations}\label{sec:simul}

\begin{figure}[t]
  \centering
  \subfigure[Gaussian]{\resizebox{2.7in}{2.4in}{
      \includegraphics{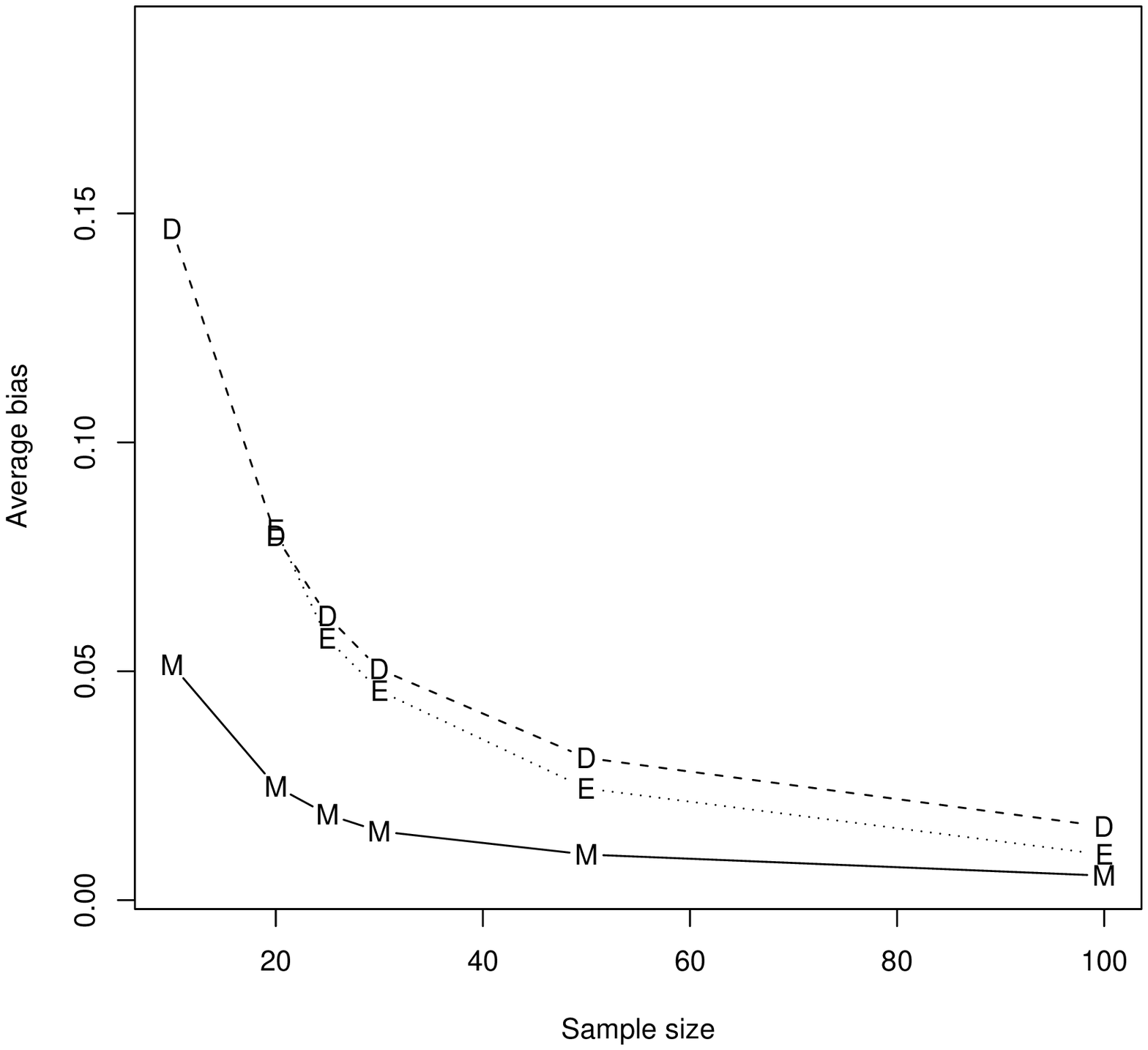}}\label{fig:bn}}
  \subfigure[$t_5$]{\resizebox{2.7in}{2.4in}{
      \includegraphics{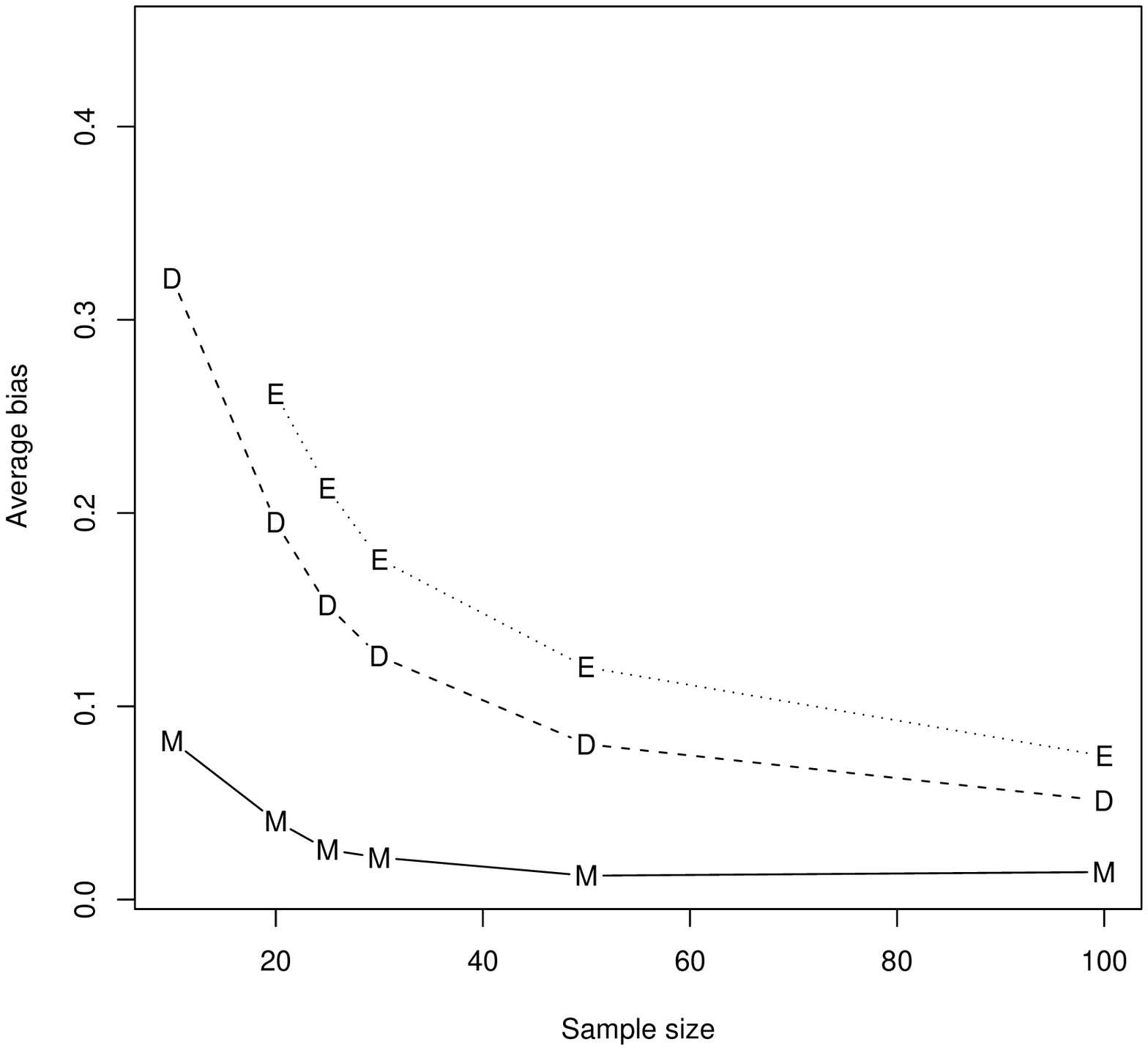}}\label{fig:b5}}
\caption{Comparing the bias of the ML, DL and the EL estimator for
various sample sizes.  ML
($\relbar\!\relbar$M$\relbar\joinrel\relbar$), DL
($\relbar\hspace{2pt}\relbar$D$\relbar\hspace{2pt}\relbar$) and EL
($\cdots$E$\cdots$).}  \label{fig:bias}
\end{figure}
\begin{figure}[t]
\centering
\subfigure[Gaussian]{
  \resizebox{2.7in}{2.4in}{
    \includegraphics{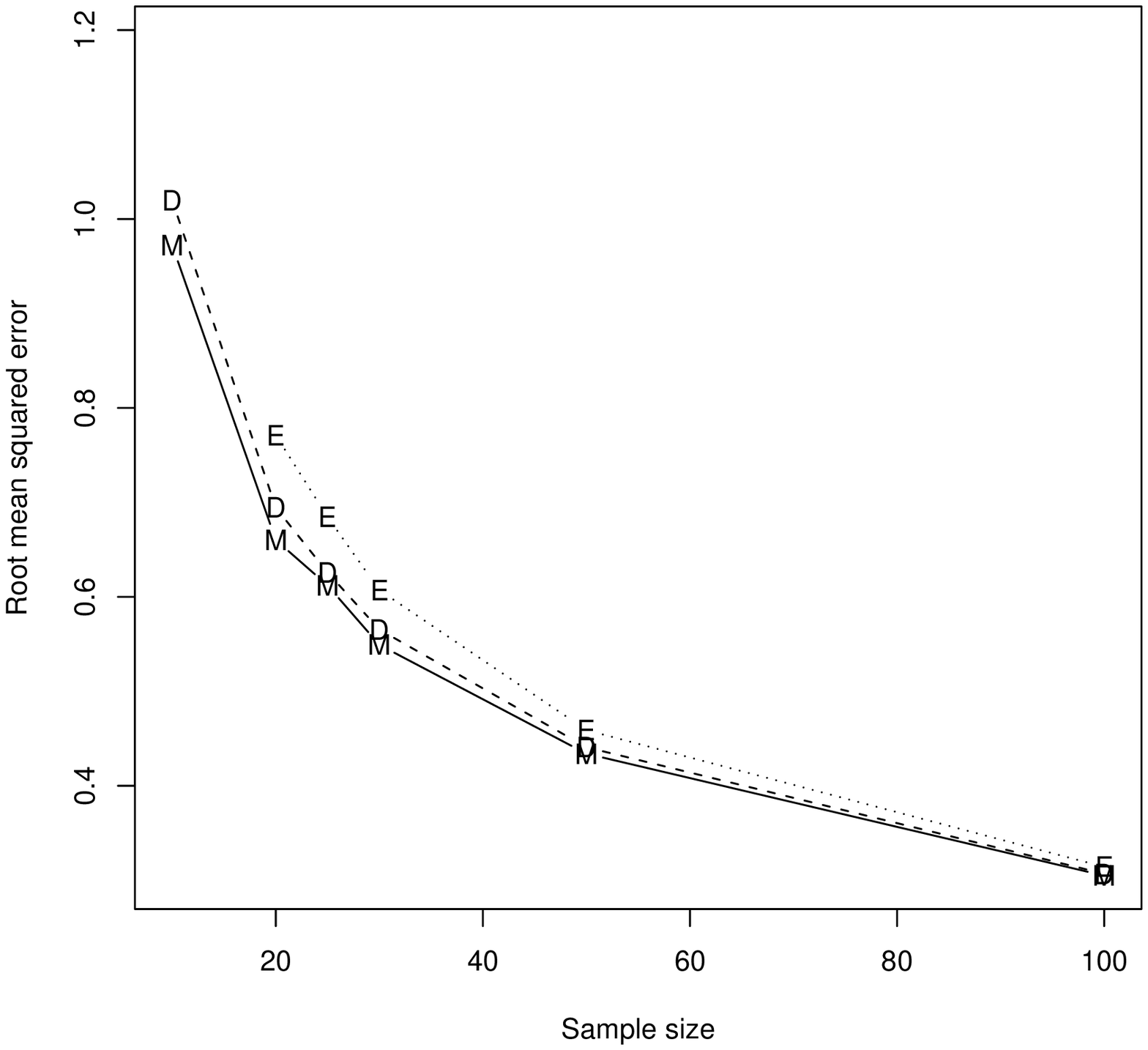}
  }
  \label{fig:l2n}}
\subfigure[$t_5$]{
  \resizebox{2.7in}{2.4in}{
    \includegraphics{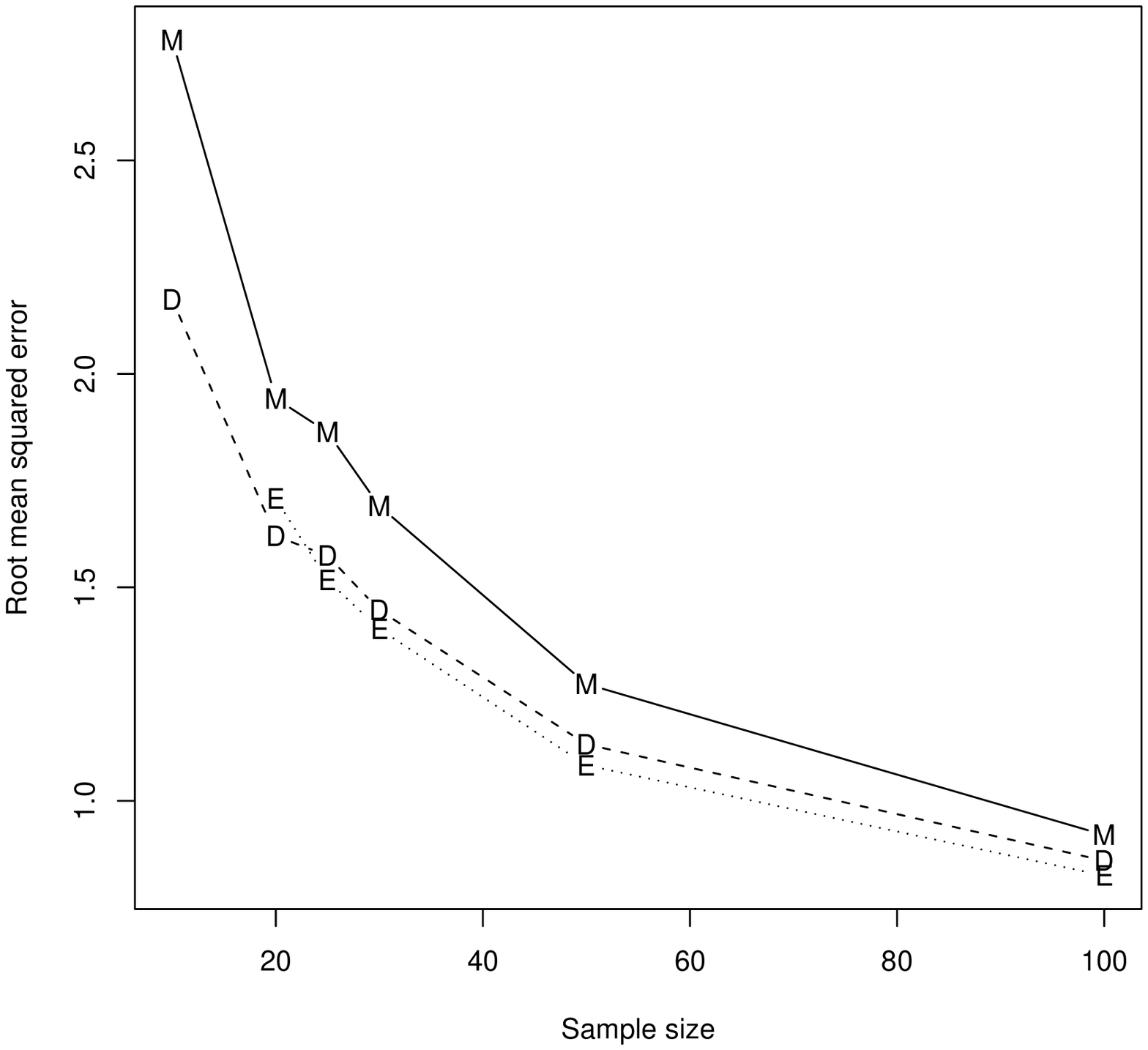}
  }
  \label{fig:l25}}
\caption{Comparing the root mean squared errors of the ML, DL and the
EL estimator for various sample sizes.  ML
($\relbar\joinrel\relbar$M$\relbar\joinrel\relbar$) DL
($\relbar\hspace{2pt}\relbar$D$\relbar\hspace{2pt}\relbar$) and
EL($\cdots$E$\cdots$).}
\label{fig:mse}
\end{figure}

Since the ML estimator $\h{\Sigma}_{M}$ and
Kauermann's dual estimator $\h{\Sigma}_D$ are based on a normality
assumption, but the empirical likelihood based estimator 
$\h{\Sigma}_{E}$ is not, it is of interest to compare their
performance, both when the underlying distribution is, and is not, Gaussian.
We simulated $1000$ data sets for sample sizes $n=10,20,25,30,50,100$ from
a multivariate normal distribution, and a multivariate $t$ distribution with
5 degrees of freedom~($t_5$).  The mean vector was zero and the covariance
matrix for the multivariate normal distribution was
\begin{equation}\label{eq:sigmasim}
\Sigma=\begin{pmatrix}
1&0&\frac{1}{2}&0\\
0&1&0&\frac{1}{4}\\
\frac{1}{2}&0&1&\frac{3}{4}\\
0&\frac{1}{4}&\frac{3}{4}&1
\end{pmatrix},
\end{equation}
corresponding to the graph shown in Figure \ref{fig:ex.bi}.  For the
$t_5$ distribution, we used $\Sigma$ as dispersion matrix, which
results in the covariance matrix $\tfrac{5}{3}\Sigma$.
In Figure \ref{fig:bias} and \ref{fig:mse} we present the bias and
root-mean-squared error (RMSE) respectively for the three estimators
(off-diagonal entries are considered once).  For the heavier-tailed
multivariate $t_5$ distribution, moments up to fourth order exist
\citep{kotz_t}, thus it makes sense to consider RMSE of the estimated
variances and covariances.  For sample size $10$ we experienced
problems with the empirical likelihood procedure, resulting from an
inability to find feasible starting values.  Consequently we do not
present results for $\h{\Sigma}_E$ when $n=10$.

From Figures \ref{fig:mse}(a) and (b) it is evident that the biases
of $\h{\Sigma}_{E}$ and $\h{\Sigma}_{D}$ are larger than the bias of
$\h{\Sigma}_{M}$ for all values of $n$.  Whereas in the
Gaussian case, $\h{\Sigma}_E$ behaves better in terms of bias than
$\h{\Sigma}_D$ for $n>20$, the opposite happens under the
$t_5$ distribution.  As would be expected the RMSE of $\h{\Sigma}_M$ is
slightly lower than 
that of $\h{\Sigma}_D$ and $\h{\Sigma}_E$ under Gaussianity;
cf.~Figure~\ref{fig:l2n}.
On the other hand $\h{\Sigma}_E$ performs better in terms of RMSE than
$\h{\Sigma}_M$ 
for all sample sizes when the underlying distribution is $t_5$; see
Figure~\ref{fig:l25}.  The RMSE of  $\h{\Sigma}_E$ 
is also smaller than that of $\h{\Sigma}_D$
under $t_5$, for moderately large sample sizes ($n>20$).



\section{Discussion}\label{extension}

We have considered three methods for estimating a covariance matrix
with pre-specified zeros.  In a Gaussian covariance graph model both ML
estimation and the dual likelihood method of \cite{kauermann:dual}
provide efficient estimates of the covariance matrix. If the
assumption of multivariate normality is not reasonable, then
non-parametric estimates can be obtained in an empirical likelihood
approach.

For the problem of maximizing the likelihood function of a Gaussian
covariance graph model we have introduced the new Iterative
Conditional Fitting (ICF) algorithm, which can be implemented in both
a univariate as well as a multivariate version.  The advantage of
multivariate ICF is the maximization of the likelihood function over
larger sections of the parameter space; the disadvantage is the
overhead in carrying out generalized least squares computations as
opposed to the standard least squares computations of univariate ICF.
Future practical experience will show whether general recommendations
in this trade off can be given, but the structure of the particular
covariance graph considered will certainly be important.

Besides its clear convergence properties, strengths of ICF include
the fact that the covariance matrix estimates are positive definite
at any stage of the algorithm and that only tools from least squares
regression are required for implementation.  In addition, it is very
appealing that ICF extends the duality between covariance graph and
undirected (concentration) graph models \cite[cf.][]{kauermann:dual}
to the level of fitting algorithms.  The commonly used method for
fitting undirected graph models, the iterative proportional fitting
(IPF) algorithm \citep[pp.182--185]{whittaker:bk}, fits marginal
distributions while fixing conditionals. ICF does exactly the
converse.  The abstract idea behind ICF can be expressed in terms of
marginal and conditional distributions which suggests that it is not
limited in any way to Gaussian covariance graph models. In fact, work
by the authors on applying ICF in binary graphical models for marginal
independence appears promising.

The ICF algorithm resembles the Iterative Conditional Modes (ICM)
algorithm of \cite{besag:1986}.  However, ICM obtains maximum {\em a
   posteriori\/} estimates in a Bayesian framework, whereas our ICF
maximizes a likelihood function, which constitutes a very differently
structured problem.  Another related algorithm is the Conditional
Iterative Proportional Fitting (CIPF) algorithm of
\cite{cramer:1998,cramer:2000}.  CIPF can be used to maximize the
likelihood function of a model that comprises joint distributions with
prescribed conditional distributions.  However, CIPF differs from ICF
because the update steps of ICF do not simply equate a conditional
distribution with a prescribed conditional, but rather maximize a
conditional likelihood function that will generally not be the same in
two different iterations of ICF.

It is obviously a most attractive feature of the empirical likelihood
procedure that it does not require multivariate normaliy.
Algorithmically, empirical likelihood estimation is more involved than
maximum likelihood and dual estimation.  In particular, we had
difficulties obtaining empirical likelihood estimates
for smaller sample sizes, which is related to a
fundamental difference between empirical likelihood estimation and the
other two methods based on multivariate normality.  Both ML and dual
estimation are possible whenever the sample covariance matrix is
positive definite, which occurs with probability one if the sample
size is larger than the number of variables, and may occur for smaller
sample sizes if the covariance graph is disconnected.  In
contrast, the optimization problem to be solved for empirical
likelihood estimation may become infeasible if the sample size is
small compared to the number of constraints imposed.  The number of
constraints depends on the covariance graph, and seemingly simpler
sparser structures impose more constraints and render the empirical
likelihood approach more sample size-demanding.

Not surprisingly, our simulations show that the ML estimates computed
with ICF are preferable if the underlying distribution is indeed
multivariate normal.  When simulating from a multivariate $t$
distribution instead non-parametric estimation via empirical
likelihood gave the best results in terms of mean squared error.

\section*{Acknowledgment}

We thank Steffen Lauritzen for pointing out the duality between ICF
and IPF, and Art Owen for suggesting use of empirical likelihood.
Sanjay Chaudhuri thanks Mark Handcock for helpful discussions.  This
work was supported by the U.S.~National Science Foundation
(DMS-9972008), the University of Washington Royalty Research Fund, the
William and Flora Hewlett Foundation, and the U.S.~National Institute of
Child Health and Human Development (R01-HD043472-01).

\bigskip\bigskip\bigskip\noindent
\sc
\parbox[t]{7cm}{
Department of Statistics \\
The University of Chicago\\
5734 S. University Avenue\\
Chicago, IL 60637\\
U.S.A.\\
E-Mail:  {\tt drton@galton.uchicago.edu}
}\hfill
\parbox[t]{7cm}{
Department of Statistics\\
University of Washington\\
Box 354322\\
Seattle, WA 98105-4322\\
U.S.A.\\
E-Mail:  {\tt \{sanjay|tsr\}@stat.washington.edu}
}

\small
\begin{thebibliography}{}

\bibitem[\protect\citeauthoryear{Anderson}{Anderson}{1969}]{anderson:1969}
Anderson, T.~W. (1969).
\newblock Statistical inference for covariance matrices with linear structure.
\newblock In {\em Multivariate Analysis, II (Proc. Second Internat. Sympos.,
  Dayton, Ohio, 1968)}, pp.\  55--66. New York: Academic Press.

\bibitem[\protect\citeauthoryear{Anderson}{Anderson}{1970}]{anderson:1970}
Anderson, T.~W. (1970).
\newblock Estimation of covariance matrices which are linear combinations or
  whose inverses are linear combinations of given matrices.
\newblock In {\em Essays in Probability and Statistics}, pp.\  1--24.
  University of North Carolina Press, Chapel Hill, N.C.

\bibitem[\protect\citeauthoryear{Anderson}{Anderson}{1973}]{anderson:1973}
Anderson, T.~W. (1973).
\newblock Asymptotically efficient estimation of covariance matrices with
  linear structure.
\newblock {\em Ann. Statist.\/}~{\em 1}, 135--141.

\bibitem[\protect\citeauthoryear{Anderson and Olkin}{Anderson and
  Olkin}{1985}]{anderson:1985}
Anderson, T.~W. and I.~Olkin (1985).
\newblock Maximum-likelihood estimation of the parameters of a multivariate
  normal distribution.
\newblock {\em Linear Algebra Appl.\/}~{\em 70}, 147--171.

\bibitem[\protect\citeauthoryear{Besag}{Besag}{1986}]{besag:1986}
Besag, J. (1986).
\newblock On the statistical analysis of dirty pictures.
\newblock {\em J. Roy. Statist. Soc. Ser. B\/}~{\em 48\/}(3), 259--302.

\bibitem[\protect\citeauthoryear{Buhl}{Buhl}{1993}]{buhl:1993}
Buhl, S. (1993).
\newblock On the existence of maximum likelihood estimators for graphical
  {G}aussian models.
\newblock {\em Scand. J. Statist.\/}~{\em 20}, 263--270.

\bibitem[\protect\citeauthoryear{Chaudhuri, Handcock, and Rendall}{Chaudhuri
  et~al.}{2005}]{scmsh1}
Chaudhuri, S., M.~S. Handcock, and M.~Rendall (2005).
\newblock Generalised linear models incorporating population level information:
  An empirical likelihood based approach.
\newblock Technical Report 484, Department of Statistics, University of
  Washington.

\bibitem[\protect\citeauthoryear{Cox and Wermuth}{Cox and
  Wermuth}{1993}]{coxwerm:lindep}
Cox, D.~R. and N.~Wermuth (1993).
\newblock Linear dependencies represented by chain graphs (with discussion).
\newblock {\em Statist. Sci.\/}~{\em 8}, 204--218,247--277.

\bibitem[\protect\citeauthoryear{Cox and Wermuth}{Cox and
  Wermuth}{1996}]{coxwerm:book}
Cox, D.~R. and N.~Wermuth (1996).
\newblock {\em Multivariate Dependencies: Models, Analysis and Interpretation}.
\newblock London: Chapman and Hall.

\bibitem[\protect\citeauthoryear{Cox, Wermuth, and Marchetti}{Cox
  et~al.}{2004}]{wermcoxmarc:04}
Cox, D.~R., N.~Wermuth, and G.~Marchetti (2004).
\newblock Decompositions and estimation of a chain of covariances.
\newblock Technical report, Department of Mathematical Statistics, Chalmers
  G{\"{o}}teborgs Universitet.

\bibitem[\protect\citeauthoryear{Cramer}{Cramer}{1998}]{cramer:1998}
Cramer, E. (1998).
\newblock Conditional iterative proportional fitting for {G}aussian
  distributions.
\newblock {\em J. Multivariate Anal.\/}~{\em 65\/}(2), 261--276.

\bibitem[\protect\citeauthoryear{Cramer}{Cramer}{2000}]{cramer:2000}
Cramer, E. (2000).
\newblock Probability measures with given marginals and conditionals:
  {$I$}-projections and conditional iterative proportional fitting.
\newblock {\em Statist. Decisions\/}~{\em 18\/}(3), 311--329.

\bibitem[\protect\citeauthoryear{Drton}{Drton}{2005}]{drton:2004b}
Drton, M. (2005).
\newblock Computing all roots of the likelihood equations of seemingly
  unrelated regressions.
\newblock {\em J. Symbolic Comput.\/}, accepted.

\bibitem[\protect\citeauthoryear{Drton and Eichler}{Drton and
  Eichler}{2005}]{drtoneichler:2004}
Drton, M. and M.~Eichler (2005).
\newblock Maximum likelihood estimation in {G}aussian chain graph models under
  the alternative {M}arkov property.
\newblock {\em Scand. J. Statist.\/}, accepted, {\tt math.ST/0508266}.

\bibitem[\protect\citeauthoryear{Drton and Perlman}{Drton and
  Perlman}{2004}]{drtonperlman:2004}
Drton, M. and M.~D. Perlman (2004).
\newblock Model selection for {G}aussian concentration graphs.
\newblock {\em Biometrika\/}~{\em 91}, 591--602.

\bibitem[\protect\citeauthoryear{Drton and Perlman}{Drton and
  Perlman}{2005}]{drtonperlman:2004b}
Drton, M. and M.~D. Perlman (2005).
\newblock A {SIN}ful approach to {G}aussian graphical model selection.
\newblock Submitted, {\tt math.ST/0508267}

\bibitem[\protect\citeauthoryear{Drton and Richardson}{Drton and
  Richardson}{2004}]{drton:2004}
Drton, M. and T.~S. Richardson (2004).
\newblock Multimodality of the likelihood in the bivariate seemingly unrelated
  regressions model.
\newblock {\em Biometrika\/}~{\em 91}, 383--392.

\bibitem[\protect\citeauthoryear{Eaton and Perlman}{Eaton and
  Perlman}{1973}]{eaton:1973}
Eaton, M.~L. and M.~D. Perlman (1973).
\newblock The non-singularity of generalized sample covariance matrices.
\newblock {\em Ann. Statist.\/}~{\em 1}, 710--717.

\bibitem[\protect\citeauthoryear{Edwards}{Edwards}{2000}]{edwards:2000}
Edwards, D.~M. (2000).
\newblock {\em Introduction to Graphical Modelling\/} (Second ed.).
\newblock New York: Springer-Verlag.

\bibitem[\protect\citeauthoryear{Gasch, Spellman, Kao, Carmel-Harel, Eisen,
  Storz, Botstein, and Brown}{Gasch et~al.}{2000}]{gasch2000}
Gasch, A.~P., P.~T. Spellman, C.~M. Kao, O.~Carmel-Harel, M.~B. Eisen,
  G.~Storz, D.~Botstein, and P.~O. Brown (2000).
\newblock Genomic expression programs in the response of yeast cells to
  environmental changes.
\newblock {\em Molecular Biology of the Cell\/}~{\em 11\/}(12), 4241--57.

\bibitem[\protect\citeauthoryear{Grzebyk, Wild, and Chouani\`ere}{Grzebyk
  et~al.}{2004}]{grzebyk:2004}
Grzebyk, M., P.~Wild, and D.~Chouani\`ere (2004).
\newblock On identification of multi-factor models with correlated residuals.
\newblock {\em Biometrika\/}~{\em 91}, 141--151.

\bibitem[\protect\citeauthoryear{Harville}{Harville}{1997}]{harville:1997}
Harville, D.~A. (1997).
\newblock {\em Matrix Algebra from a Statistician's Perspective}.
\newblock New York: Springer-Verlag.

\bibitem[\protect\citeauthoryear{Hellerstein and Imbens}{Hellerstein and
  Imbens}{1999}]{hellimbens}
Hellerstein, J. and G.~W. Imbens (1999).
\newblock Imposing moment restrictions from auxiliary data by weighting.
\newblock {\em The Review of Economics and Statistics\/}~{\em LXXXI\/}(1),
  1--14.

\bibitem[\protect\citeauthoryear{Kauermann}{Kauermann}{1996}]{kauermann:dual}
Kauermann, G. (1996).
\newblock On a dualization of graphical {G}aussian models.
\newblock {\em Scand. J. Statist.\/}~{\em 23}, 105--116.

\bibitem[\protect\citeauthoryear{Kotz and Nadarajah}{Kotz and
  Nadarajah}{2004}]{kotz_t}
Kotz, S. and S.~Nadarajah (2004).
\newblock {\em Multivariate {$t$} Distributions and their Applications}.
\newblock Cambridge: Cambridge University Press.

\bibitem[\protect\citeauthoryear{Lauritzen}{Lauritzen}{1996}]{lau:bk}
Lauritzen, S.~L. (1996).
\newblock {\em Graphical Models}.
\newblock Oxford, UK: Clarendon Press.

\bibitem[\protect\citeauthoryear{Mao, Kschischang, and Frey}{Mao
  et~al.}{2004}]{mao:2004}
Mao, Y., F.~R. Kschischang, and B.~J. Frey (2004).
\newblock Convolutional factor graphs as probabilistic models.
\newblock In U.~Kj{\ae}rulff and C.~Meek (Eds.), {\em Procceding of the 20th
  Conference on Uncertainty in Artificial Intelligence}, pp.\  374--381. San
  Francisco: Morgan Kaufmann.

\bibitem[\protect\citeauthoryear{Owen}{Owen}{2001}]{owen:2001}
Owen, A.~B. (2001).
\newblock {\em Empirical Likelihood}.
\newblock Boca Raton: Chapman \& Hall.

\bibitem[\protect\citeauthoryear{Qin and Lawless}{Qin and
  Lawless}{1994}]{Qinlawless1}
Qin, J. and J.~Lawless (1994).
\newblock Empirical likelihood and general estimating equations.
\newblock {\em Ann. Statist.\/}~{\em 22}, 300--325.

\bibitem[\protect\citeauthoryear{Richardson and Spirtes}{Richardson and
  Spirtes}{2002}]{richardson:2002}
Richardson, T.~S. and P.~Spirtes (2002).
\newblock Ancestral graph {M}arkov models.
\newblock {\em Ann.~Statist.\/}~{\em 30}, 962--1030.

\bibitem[\protect\citeauthoryear{Whittaker}{Whittaker}{1990}]{whittaker:bk}
Whittaker, J. (1990).
\newblock {\em Graphical Models in Applied Multivariate Statistics}.
\newblock Chichester: Wiley.

\bibitem[\protect\citeauthoryear{Wright}{Wright}{1921}]{wright:cnc}
Wright, S. (1921).
\newblock Correlation and {C}ausation.
\newblock {\em J. Agricultural Research\/}~{\em 20}, 557--585.

\bibitem[\protect\citeauthoryear{Zellner}{Zellner}{1962}]{zellner:62}
Zellner, A. (1962).
\newblock An efficient method of estimating seemingly unrelated regression
  equations and tests for aggregation bias.
\newblock {\em J. Amer. Statist. Assoc.\/}~{\em 57}, 348--368.

\end{thebibliography}
\end{document}